\documentclass[11pt]{article}
\usepackage{ graphics}

\newcommand{\mathsym}[1]{{}}
\newcommand{\unicode}[1]{{}}
\usepackage{latexsym,xy,eucal,mathrsfs,graphs}
\usepackage{setspace}
\usepackage{amsthm}
\usepackage{pstricks}
\usepackage{amssymb}
\usepackage{amsmath}
\usepackage{indentfirst}
\usepackage{epsfig}
\usepackage{indentfirst}
\usepackage{amscd}
\usepackage{multirow}

\newtheorem{theorem}{Theorem}[section]

\newtheorem{lemma}[theorem]{Lemma}
\newtheorem{claim}[theorem]{CLAIM}

\newtheorem{observation}[theorem]{Observation}
\newtheorem{result}[theorem]{Result}
\theoremstyle{definition}
\newtheorem{defn}[theorem]{Definition}
\newtheorem{example}[theorem]{Example}

\usepackage{simplemargins}
\setleftmargin{0.88in}
\setrightmargin{0.88in}
\settopmargin{0.817in}
\setbottommargin{0.83in}

 \font\conj =
msbm10 at 12pt  

\def\text#1{\mathop{{\rm{#1}}}\nolimits}
\def\mr{\mathop{{\rm{mr}}}\nolimits}
\def\msr{\mathop{{\rm{msr}}}\nolimits}
\def\max{\mathop{{\rm{max}}}\nolimits}

\def\nulli{\mathop{{\rm{null}}}\nolimits}

\def\cc{\mathop{{\rm{cc}}}\nolimits}

\def\d{\mathop{{\rm{d}}}\nolimits}
\def\sp{\mathop{{\rm{Span}}}\nolimits}

\def\Re{{\mbox{\conj R}}}

\def\F{{\mbox{\conj F}}}
\def\C{{\mbox{\conj C}}}
\def\Q{{\mbox{\conj Q}}}
\def\dis{\displaystyle}
\def\1{\'{\i}}

\title{ On $\delta$-Graphs and  Delta Conjecture}

\author{Pedro D\'{\i}az Navarro\thanks{Escuela de Matem\'atica, Universidad de Costa Rica}}
\date{August, 2016}
\begin{document}
\maketitle

\begin{abstract}
\noindent
In  this paper   we  define  two infinite families  of graphs  called C-$\delta$ graphs  and  $\delta$- graph  and  prove  that   $\delta$-graphs  satisfy  $\delta$  conjecture. Also   we  introduce   a  family  of  C-$\delta$ graphs  from  which   we can  identify   $\delta$  graphs  as  their  complements.  Finally   we   give a list of  C-$\delta$ graphs   and  the  relationship  with  their  minimum  semidefinite  rank.
\end{abstract}

\section{Introduction}

A {\it graph} $G$ consists of a set  of vertices $V(G)=\{1,2,\dots,n\}$ and a set of edges $E(G)$, where an edge is defined to be an unordered  pair of vertices. The {\it order} of $G$, denoted $\vert G\vert $,  is the cardinality  of $V(G)$. A graph is {\it simple} if it has no multiple  edges  or loops. The {\it complement } of a graph $G(V,E)$ is the graph $\overline{G}=(V,\overline{E})$, where $\overline{E}$ consists of all those edges of the complete  graph $K_{\vert G\vert}$ that are not in $E$.

 A matrix $A=[a_{ij}]$ is {\it combinatorially symmetric} when $a_{ij}=0$ if and only if $a_{ji}=0$. We say that  $G(A)$  is the graph of a combinatorially symmetric matrix $A=[a_{ij}]$ if $V=\{1,2,\dots,n\}$ and $E=\{\{i,j\}: a_{ij}\ne0\}$ . The main diagonal entries of $A$ play no role in determining $G$. Define $S(G,\F)$ as the set of all $n\times n$ matrices that are real symmetric if $\F=\Re$ or complex Hermitian if $\F=\C$ whose graph is $G$. The sets $S_+(G,\F)$ are the corresponding subsets of positive semidefinite (psd) matrices. The smallest possible rank  of any matrix $A\in S(G,\F)$  is the {\it minimum rank} of $G$, denoted $\mr(G,\F)$, and the smallest possible rank of any matrix $A\in S_+(G,\F)$  is the {\it minimum semidefinite rank} of $G$, denoted $\mr_+(G)$ or $\msr(G)$.

In 1996, the minimum rank among real symmetric matrices with a given graph was studied  by  Nylen \cite{PN}. It gave rise to the area of minimum rank problems which led to the study of minimum rank among complex Hermitian matrices and positive semidefinite matrices associated with a given graph. Many results can be  found for  example  in  \cite{FW2, VH, YL, LM, PN}.

During  the   AIM workshop of 2006 in Palo Alto, CA, it  was  conjectured   that for any   graph $G$  and  infinite field  $F$, $\mr(G,\F)\le |G|-\delta(G)$ where $\delta(G)$  is the minimum  degree  of  $G$.  It  was  shown  that  for if  $\delta(G)\le 3$ or  $\delta(G)\ge |G|-2$    this inequality  holds.  Also  it  can be verified that  if  $|G|\le 6$ then  $\mr(G,F)\le |G|-\delta(G)$.  Also  it  was   proven  that any  bipartite graph  satisfies  this conjecture. This  conjecture   is called  the {\it Delta Conjecture}. If we   restrict  the study  to  consider  matrices  in $S_+(G,\F)$  then  delta conjecture  is  written  as $\msr(G)\le |G|-\delta(G)$. Some  results on delta conjecture  can be found in \cite{AB,RB1, SY1, SY} but  the   general problem  remains  unsolved.

\section{Graph Theory Preliminaries}
\addtocontents{toc}{\vspace{-15pt}}
In  this section   we give   definitions and results from  graph theory which    will  be used in  the remaining  chapters. Further details  can be found  in \cite{BO,BM, CH}.

A {\bf graph} {$G(V,E)$} is a pair {$(V(G),E(G)),$} where {$V(G)$} is the set of vertices and {$E(G)$} is the set of edges together  with an  {\bf  incidence  function} $\psi(G)$ that associate with  each edge  of  $G$ an  unordered  pair  of (not necessarily  distinct) vertices  of  $G$. The {\bf order} of {$G$}, denoted {$|G|$}, is the number of vertices in {$G.$} A graph is said to be {\bf simple} if it has no loops or multiple edges. The {\bf complement} of a graph {$G(V,E)$} is the graph {$\overline{G}=(V,\overline{E}),$} where {$\overline{E}$} consists of all the edges that are not in {$E$}.
A {\bf  subgraph} {$H=(V(H),E(H))$} of {$G=(V,E)$} is a graph with {$V(H)\subseteq V(G)$} and {$E(H)\subseteq E(G).$} An {\bf induced subgraph} {$H$} of {$G$}, denoted G[V(H)], is a subgraph with {$V(H)\subseteq V(G)$} and {$E(H)=\{\{i,j\} \in E(G):i,j\in V(H)\}$}. Sometimes  we  denote the  edge $\{i,j\}$ as $ij$.

We  say  that  two  vertices of a graph $G$  are {\bf adjacent}, denoted  $v_i\sim v_j$,   if  there is an edge $\{v_i,v_j\}$ in  $G$.  Otherwise  we say  that the  two  vertices $v_i$ and $v_j$ are {\bf non-adjacent}  and  we denote this  by $v_i \not\sim v_j$.  Let {$N(v)$} denote the set of vertices that are adjacent to the vertex {$v$} and let {$N[v]=\{v\}\cup N(v)$}. The {\bf degree} of a vertex {$v$} in {$G,$} denoted {$\d_G(v),$} is the cardinality of {$N(v).$} If {$\d_G(v)=1,$} then {$v$} is said to be a {\bf pendant} vertex of {$G.$} We use {$\delta(G)$} to denote the minimum degree of the vertices in {$G$}, whereas {$\Delta(G)$} will denote the maximum degree of the vertices in {$G$}.
 
 Two  graphs $G(V,E)$ and $H(V',E')$  are  identical  denoted  $G=H$, if  $V=V',  E=E'$, and $\psi_G=\psi_H$  . Two  graphs $G(V,E)$ and $H(V',E')$ are {\bf isomorphic}, denoted  by  $G\cong H$, if  there exist bijections $\theta:V\to V'$  and $\phi: E\to  E' $ such  that $\psi_G(e)=\{u,v\}$  if  and  only if  $\psi_H(\phi(e))= \{\theta(u), \theta(v)\}$.
 
 A {\bf complete graph}  is a simple graph in which the vertices are pairwise adjacent.
  We will use {$nG$} to denote {$n$} copies of a graph {$G$}. For example, $3K_1$  denotes three  isolated vertices $K_1$ while {$2K_2$} is the graph given  by   two  disconnected  copies  of $K_2$.
  
 A {\bf path} is a list of distinct vertices in which successive vertices are connected by edges. A path on {$n$} vertices is denoted by {$P_n.$} A graph {$G$} is said to be {\bf connected} if there is a path between any two vertices of {$G$}. A {\bf cycle} on {$n$} vertices, denoted {$C_n,$} is a path such that the beginning vertex and the end vertex are the same. A {\bf tree} is a connected graph with no cycles. A graph $G(V,E)$ is  said  to be {\bf chordal}  if it has no induced cycles $C_n$ with $n\ge 4$.
 A  {\bf component}  of a graph $G(V,E)$ is  a maximal connected  subgraph. A  {\bf cut vertex}   is  a vertex  whose deletion  increases  the number of  components.
 
 The {\bf union}  $G\cup G_2$ of  two  graphs $G_1(V_1,E_1)$ and $G_2(V_2,G_2)$  is  the union  of  their  vertex  set  and  edge  set,  that is $G\cup G_2(V_1\cup V_2,E_1\cup E_2$. When  $V_1$ and $V_2$ are disjoint their union  is called  {\bf  disjoint union} and  denoted $G_1\sqcup G_2$.

%******************************************************************************
%                       ORTHOGONAL REPRESENTATION
%******************************************************************************

\section{The Minimum  Semidefinite Rank  of  a Graph}
%\addtocontents{toc}{\vspace{-15pt}}
In  this section  we will establish   some of  the   results for  the minimum  semidefinite rank ($\msr$)of a graph $G$  that  we  will be using in the subsequent chapters.

A {\bf positive  definite} matrix  $A$ is an Hermitian   $n\times n$ matrix such  that $x^\star A x>0$  for all nonzero  $x\in \C^n$. Equivalently,  $A$  is a  $n\times n$ Hermitian positive definite matrix  if and  only  if   all the  eigenvalues of $A$ are positive (\cite{RC}, p.250).

A $n\times n$ Hermitian matrix  $A$ such  that $x^\star A x\ge 0$  for all $x\in \C^n$ is  said  to be  {\bf positive  semidefinite (psd)}. Equi\-va\-lently,   $A$ is a  $n\times n$ Hemitian positive  semidefinite matrix if and  only  if  $A$  has  all   eigenvalues nonnegative (\cite{RC}, p.182).

If $\overrightarrow{V}=\{\overrightarrow{v_1},\overrightarrow{v_2},\dots, \overrightarrow{v_n}\}\subset \Re^m$  is a set of  column vectors  then  the  matrix
$ A^T A$, where $A= \left[\begin{array}{cccc}
  \overrightarrow{v_1} & \overrightarrow{v_2} &\dots& \overrightarrow{v_n}
\end{array}\right]$
and $A^T$  represents  the  transpose matrix of  $A$, is a psd matrix  called  the  {\bf Gram matrix} of $\overrightarrow{V}$. Let $G(V,E)$  be a graph  associated  with  this Gram matrix. Then  $V_G=\{v_1,\dots, v_n\}$ correspond to  the set of  vectors in $\overrightarrow{V}$ and  E(G) correspond to  the nonzero inner products  among  the  vectors  in $\overrightarrow{V}$. In this  case $\overrightarrow{V}$  is  called an  {\bf orthogonal representation} of $G(V,E)$ in $\Re^m$. If  such  an  orthogonal  representation  exists  for  $G$ then $\msr(G)\le m$.

The  {\bf  maximum positive  semidefinite   nullity }  of a graph $G$, denoted $M_+(G)$   is   defined  by
$
M_+(G) =\max\{\nulli(A): A\ \hbox{ is symmetric and positive semidefinite and }\   G(A) =G \}
$,
where $G(A)$  is  the graph obtained  from  the matrix $A$. From  the  rank-nullity theorem we get  $ \msr(G)+ M_+(G)=|G|$.

Some  of  the most  common   results about  the minimum semidefinite  rank of a graph  are  the  following:

\begin{result}\cite{VH}\label{msrtree}
If  $T $ is a tree  then $\msr(T)= |T|-1$.
\end{result}
\begin{result}\cite{MP3}\label{msrcycle}
%p25
 The cycle $C_n$ has minimum semidefinite rank $n-2$.
\end{result}

% from Mathew Booth et all .On the minimum  semidefinite rank of a simple graph

% From Mathew booth et all On the minimum  rank  among  positive semidefinite  matrices with  a given  graph
\begin{result}\label{res2}
 \cite{MP3} \ If a connected  graph $G$  has a pendant  vertex $v$, then $\msr(G)=\msr(G-v)+1$ where $G-v$  is obtained as an induced subgraph of $G$ by  deleting $v$.
\end{result}

\begin{result} \cite{PB} \label{OS2}
 If {$G$} is a connected, chordal graph, then {$\msr(G)=\cc(G).$}
\end{result}
%***************************************************************************

\begin{result}\label{res1}
\cite{MP2}\ If a graph $G(V,E)$ has a cut vertex, so that $G=G_1\cdot G_2$, then  $\msr(G)= \msr(G_1)+\msr(G_2)$.
\end{result}
\section{$\mathbf \delta$-Graphs and  the Delta Conjecture}

In  this  section  we  define  a new  family of  graphs  called $\delta$-graphs  and  show  that  they satisfy  the  delta conjecture.

\begin{defn}\label{ccpg}
Suppose  that $G=(V,E)$  with $|G|=n \ge 4$  is simple  and  connected such  that  $\overline{G}=(V,\overline{E})$  is also  simple  and  connected. We  say  that  $G$ is a {\bf $\mathbf{\delta}$-graph} if  we  can  label  the vertices of $G$ in such a way that
\begin{enumerate}
  \item[(1)] the  induced    graph    of  the  vertices  $v_1,v_2,v_3$  in $G$ is  either $3K_1$  or  $K_2 \sqcup K_1$,  and
  \item[(2)] for $m\ge 4$,  the vertex  $v_m$  is  adjacent   to all   the  prior  vertices $v_1,v_2,\dots,v_{m-1}$  except  for at most $\dis{\left\lfloor\frac{m}{2}-1\right\rfloor}$ vertices.
  \end{enumerate}
 \end{defn}
\begin{defn} Suppose   that  a graph $G(V,E)$  with $|G|=n \ge 4$  is  simple and  connected   such  that  $\overline{G}=(V,\overline{E})$  is also   simple and  connected.  We  say  that  $G(V,E)$  is a {\bf C-$\mathbf{\delta}$  graph}  if    $\overline{G}$   is a  $\delta$-graph.

In other  words,  $G$   is a  {\bf C-$\mathbf{\delta}$ graph} if  we can  label   the  vertices   of  $G$  in  such a  way   that
 \begin{enumerate}
 \item[ (1)] the  induced  graph  of  the vertices  $v_1,v_2,v_3$ in  $G$ is  either  $K_3$ or $P_3$,  and
 \item[(2)]  for  $m\ge 4$,  the vertex $v_m$  is adjacent  to at most  $\dis{\left\lfloor\frac{m}{2}-1\right\rfloor}$ of  the  prior  vertices $v_1,v_2,\dots,v_{m-1}$.
\end{enumerate}
\end{defn}
\begin{example}\label{examplecp}
The   cycle $C_n, n\ge 6$  is a C-$\delta$ graph and  its  complement  is  a  $\delta$-graph.
\begin{center}
\includegraphics[height=25mm]{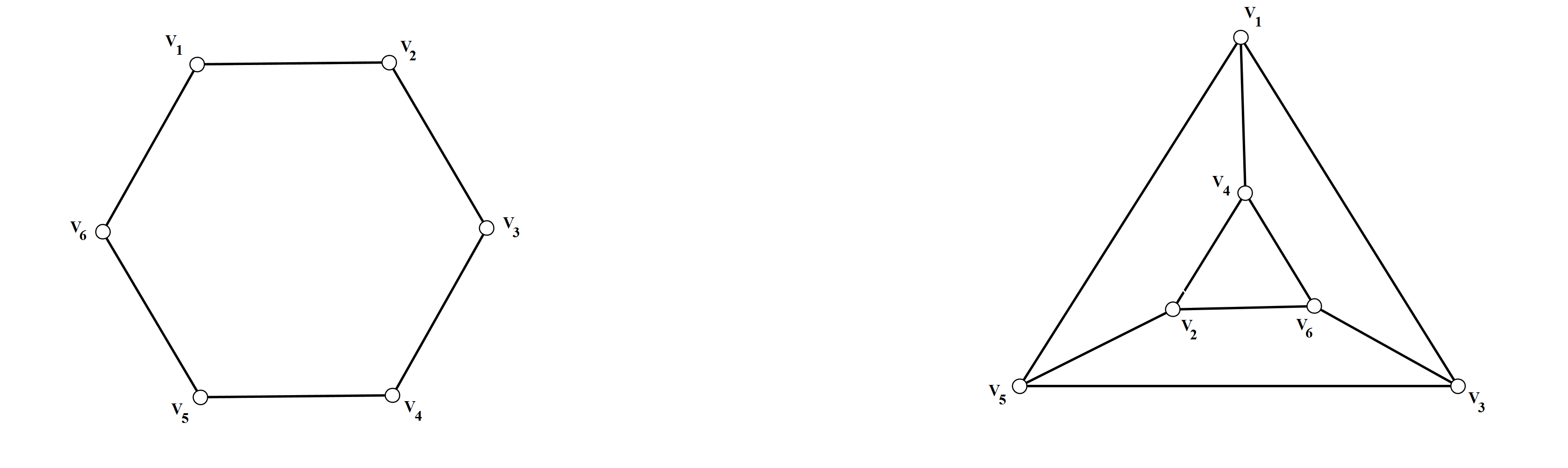}
 \end{center}
\vspace{-0.1in}\begin{figure}[h]
\centering
\caption{The  Graphs $C_6$  and  the  3-Prism }\label{fig3.1}
\end{figure}
%\vspace*{-0.9\baselineskip}
Note that we  can label  the vertices  of $C_6$ clockwise $v_1,v_2,v_3,v_4,v_5,v_6$. The  graph  induced   by  $v_1,v_2,v_3$  is $P_3$. The    vertex $v_4$ is   adjacent to a prior vertex   which  is  $v_3$. Also,   the  vertex $v_5$ is adjacent to  vertex $v_4$ and   the  vertex $v_6$ is adjacent to two prior  vertices  $v_1$ and $v_5$.  Hence, $C_6$  is C-$\delta$ graph. The $3$-prism which  is  isomorphic to  the complement  of  $C_6$,  is a $\delta$-graph.
\end{example}
\begin{lemma}\label{lem2} Let {$G(V,E)$} be a $\delta$-graph. Then   the induced  graph of {$\{v_1,v_2,v_3\}$} in {$G$} denoted  by $H$ has an  orthogonal  representation  in {$\Re^{\Delta(\overline{G})+1}$} satisfying the  following  conditions:
\begin{enumerate}
\item [(i)] the  vectors   in  the orthogonal  representation  of  $H$ can be chosen  with nonzero coordinates, and
\item [(ii)]\label{L1}$\overrightarrow{v}\not\in \sp(\overrightarrow{u})$ for each pair   of  distinct  vertices  $u,v$ in  $H$.
\end{enumerate}
\end{lemma}
\noindent\MakeUppercase{Proof.}Let {$G(V,E)$} be a $\delta$-graph. Label  the  vertices  of  $G$ in  such a way  that the labeling  satisfies  the conditions (1) and (2) for  $\delta$-graphs. Let {$H$}   be the  induced  graph   in $G$  of  $\{v_1,v_2,v_3\}\subseteq V$. Then   $H$ is  either  $3K_1$ or $K_2\sqcup K_1$.
Since  {$G$}  and {$\overline{G}$} are simple and connected  it  follows  that
\begin{center}\label{ine1}
$2\le \Delta(\overline{G})\le n-2$
\end{center}
Let $\{\overrightarrow{e}_j\}, j=1,2,\dots,\Delta(\overline{G})+1$ be  the standard  orthonormal  basis   for $\Re^{\Delta(\overline{G})+1}$.\\
\noindent{\it Case 1.} Suppose   the induced  graph $H$ of {$\{v_1,v_2,v_3\}\subseteq V$} in $G$  is $3K_1$  which is  disconnected. Choose  $\overrightarrow{v_1},\overrightarrow{v_2},\overrightarrow{v_3}$  in  $\Re^{\Delta(\overline{G})+1}$ corresponding to  {$v_1,v_2,v_3$} respectively such  that:
\begin{eqnarray*}
% \nonumber to remove numbering (before each equation)
  \overrightarrow{v_1} &=& \dis{\sum_{j=1}^{\Delta(\overline{G})+1}k_{1,j}\overrightarrow{e}_j}\\
 \overrightarrow{v_2} &=&  \dis{\sum_{j=1}^{\Delta(\overline{G})+1}k_{2,j}\overrightarrow{e}_j} \\
 \overrightarrow{v_3} &=&  \dis{\sum_{j=1}^{\Delta(\overline{G})+1}k_{3,j}\overrightarrow{e}_j}
\end{eqnarray*}
\noindent where  the  scalars $k_{1,j}, j=1,2,\dots,\Delta(\overline{G})+1$, and $k_{2,s}, s=1,2,\dots, \Delta(\overline{G})$   are  chosen  not  zero from   different  field  extensions in the following  way:
  \begin{itemize}
 \item [] $k_{1,1}\not \in \Q$,
  \item[] $k_{1,2}\not \in\Q[k_{1,1}]$,
  \item [] $k_{1,3}\not \in \Q[k_{1,1}, k_{1,2}]$,
  \item[] $\vdots$
  \item [] $ k_{1,\Delta(\overline{G})+1}\not \in  \Q[k_{1,1},k_{1,2},\dots,k_{1 ,\Delta(\overline{G})}]$,
 \item []$ k_{2,1}\not \in \Q[k_{1,1},k_{1,2},\dots,k_{1,\Delta(\overline{G})+1}]$,
 \item[] $k_{2,2}\not \in \Q[k_{1,1},k_{1,2},\dots,k_{1,\Delta(\overline{G})+1},k_{2,1}]$,
 \item[] $\vdots$
 \item [] $k_{2,\Delta(\overline{G})}\not \in \Q[k_{1,1},k_{1,2},\dots,k_{1,\Delta(\overline{G})+1},k_{2,1},\dots,k_{2,\Delta(\overline{G})-1}]$.
 \end{itemize}
Now choose
 $$
 k_{2,\Delta(\overline{G})+1}=\frac{-1}{k_{1,\Delta(\overline{G})+1}}\dis{\sum_{j=1}^{\Delta(\overline{G})}k_{1,j}k_{2,j}}
 $$
 As  a consequence $\langle \overrightarrow{v}_1,\overrightarrow{v}_2\rangle=0$.

In  order  to  find a vector  $\overrightarrow{v}_3$,  we need to  solve  the $ 2\times (\Delta(\overline{G})+1)$ system satisfying
\begin{eqnarray*}
% \nonumber to remove numbering (before each equation)
  \langle\overrightarrow{v}_1,\overrightarrow{v}_3\rangle &=& 0 \\
  \langle\overrightarrow{v}_2,\overrightarrow{v}_3\rangle &=& 0
\end{eqnarray*}
\noindent in the variables $k_{3,j}, j=1,2,\dots,\Delta(\overline{G})+1$. The homogeneous system  has  infinitely many solutions  because $\Delta(\overline{G})+1\ge3$. Reducing the matrix of  this system to  echelon  form we  get
$$
\left(
  \begin{array}{ccccccc}
    k_{1,1} & k_{1,2}& k_{1,3}& \dots&  k_{1,\Delta(\overline{G})+1}  \\
    k_{2,1} & k_{2,2}& k_{2,3}& \dots&  k_{2,\Delta(\overline{G})+1} \\
  \end{array}
\right)\sim
$$
$$
\left(
  \begin{array}{ccccccc}
    1 & \frac{k_{1,2}}{k_{1,1}}&\frac{ k_{1,3}}{k_{1,1}}& \dots&  \frac{k_{1,\Delta(\overline{G})+1}}{k_{1,1}}  \\
    k_{2,1} & k_{2,2}& k_{2,3}& \dots&  k_{2,\Delta(\overline{G})+1} \\
  \end{array}
\right)\sim
$$

$$
\left(
  \begin{array}{ccccccc}
    1 & \frac{k_{1,2}}{k_{1,1}}&\frac{ k_{1,3}}{k_{1,1}}& \dots&  \frac{k_{1,\Delta(\overline{G})+1}}{k_{1,1}}  \\
   0 & k_{2,2}-\frac{k_{2,1}k_{1,2}}{k_{1,1}}& k_{2,3}-\frac{k_{2,1}k_{1,3}}{k_{1,1}}& \dots&  k_{2,\Delta(\overline{G})+1}-\frac{k_{2,1}k_{1,\Delta(\overline{G})+1}}{k_{1,1}} \\
  \end{array}
\right)
$$
Let $\alpha = k_{2,2}-\frac{k_{2,1}k_{1,2}}{k_{1,1}}$.  Since $\alpha \ne 0$  because  $k_{2,2}\not\in \Q[k_{1,1},\dots,k_{1,\Delta(\overline{G})+1},k_{2,1}]$. We can  get  the echelon  form  of  the matrix by multiplying  the  second  row  by  $\frac{1}{\alpha}$.
$$
\left(
  \begin{array}{ccccccc}
    1 & \frac{k_{1,2}}{k_{1,1}}&\frac{ k_{1,3}}{k_{1,1}}& \dots&  \frac{k_{1,\Delta(\overline{G})+1}}{k_{1,1}}  \\
   0 & 1&\frac{1}{ \alpha}(k_{2,3}-\frac{k_{2,1}k_{1,3}}{ k_{1,1}})& \dots&\frac{1}{ \alpha}  (k_{2,\Delta(\overline{G})+1}-\frac{k_{2,1}k_{1,\Delta(\overline{G})+1}}{ k_{1,1}}) \\
  \end{array}
\right)
$$
Since  the system has infinitely many solutions, $k_{3,j}, j= 3,\dots,\Delta(\overline{G})+1$ are free  parameters.
We can choose them from different  field  extensions  in the following  way,
\begin{itemize}
  \item[] $ k_{3,3}\not \in \Q[k_{1,1},k_{1,2},\dots,k_{1,\Delta(\overline{G})+1},k_{2,1},\dots,k_{2,\Delta(\overline{G})+1}]$,
    \item[] $ k_{3,4}\not \in \Q[k_{1,1},k_{1,2},\dots,k_{1,\Delta(\overline{G})+1},k_{2,1},\dots,k_{2,\Delta(\overline{G})+1},k_{3,1}]$,
  \item[]$ \vdots$
 \item[] $ k_{3,\Delta(\overline{G})+1}\not \in \Q[k_{1,1},k_{1,2},\dots,k_{1,\Delta(\overline{G})+1},k_{2,1},\dots,k_{2,\Delta(\overline{G})+1},k_{3,3},\dots,k_{3,\Delta(\overline{G})}]$.
  \end{itemize}
Since $k_{2,j}-\frac{k_{2,1}k_{1,j}}{ k_{1,1}}\ne 0, j=3,\dots,\Delta(\overline{G})+1$,  we can choose these parameters   such that $k_{3,1}$  and $k_{3,2}$ are also nonzero.
Therefore we get  $\langle\overrightarrow{v}_1,\overrightarrow{v}_3\rangle=\langle\overrightarrow{v}_2,\overrightarrow{v}_3\rangle=0$.
As  a  result $\overrightarrow{v_1},\overrightarrow{v_2},\overrightarrow{v_3}$ is an orthogonal representation  of the induced  graph  $H=3K_1$ satisfying conditions (i) and (ii).\\
\noindent{\it Case 2.} Suppose the induced  graph $H$ is $K_2\sqcup K_1$. Let us  assume $v_1 \not\sim v_2,v_2 \not\sim v_3$, and $v_1  \sim v_3 $.\\
The  vectors $\overrightarrow{v}_1$ and  $\overrightarrow{v}_2$  for   the  vertices  $v_1$ and $v_2$  respectively can be chosen  in the same  way as in  case 1  to  get $\langle\overrightarrow{v}_1,\overrightarrow{v}_2\rangle=0$.
In order to find   a vector $\overrightarrow{v_3}=k_{3,1}\overrightarrow{e_1}+k_{3,2}\overrightarrow{e_2}+\dots +k_{3,\Delta(\overline{G})+1}\overrightarrow{e}_{\Delta(\overline{G})+1}\in \Re^{\Delta(\overline{G})+1}$  with  nonzero components for  the vertex  $v_3$ we know  that  the vector  $\overrightarrow{v}_3$  should  satisfy the system
 \begin{eqnarray*}
 % \nonumber to remove numbering (before each equation)
 \langle\overrightarrow{v_1},\overrightarrow{v_3}\rangle&=&g_1 ,  \ \ g_1\ne 0 \label{s11} \\
  \langle\overrightarrow{v_2},\overrightarrow{v_3}\rangle &=& 0     \label {s12}
 \end{eqnarray*}
\noindent in  the variables $k_{3,j}, j=1,2,\dots,\Delta(\overline{G})+1$ because  $v_1\sim v_3$ and $v_2\not\sim v_3$ in $G$.
Therefore,  rewriting  the system   in  the form
\begin{align}
   % \nonumber to remove numbering (before each equation)
    k_{1,1}k_{3,1}+k_{1,2}k_{3,2}+\dots+k_{1,\Delta(\overline{G})+1}k_{3,\Delta(\overline{G})+1}&= g_1, \label{esys1}\\
     k_{2,1}k_{3,1}+k_{2,2}k_{3,2}+\dots+k_{2,\Delta(\overline{G})+1}k_{3,\Delta(\overline{G})+1} &= 0 \label{esys2}
   \end{align}
\noindent the augmented matrix  of the non-homogeneous system  becomes
$$
\left(
  \begin{array}{ccccccc}
    k_{1,1} & k_{1,2}& k_{1,3}& \dots&  k_{1,\Delta(\overline{G})+1} & | & g_1 \label{ms1} \\
    k_{2,1} & k_{2,2}& k_{2,3}& \dots&  k_{2,\Delta(\overline{G})+1} & |  & 0 \label {ms2}\\
  \end{array}
\right)
$$
Since $\Delta(\overline{G})+1\ge 3$,  the system   has infinitely many  solutions  depending  on  one or  more  free parameters if  the system is consistent. Since  $k_{1,1}$  is nonzero,  we can  divide  the first   row by $k_{1,1}$ to obtain
$$
\left(
  \begin{array}{ccccccc}
   1 &\frac{ k_{1,2}}{k_{1,1}}&\frac{k_{1,3}}{k_{1,1}}& \dots&  \frac{k_{1,\Delta(\overline{G})+1} }{k_{1,1}}& | & \frac{g_1}{k_{1,1}} \label{ms3} \\
    k_{2,1} & k_{2,2}& k_{2,3}& \dots&  k_{2,\Delta(\overline{G})+1} & |  & 0 \label {ms4}\\
  \end{array}
\right)\sim
$$
$$
\left(
  \begin{array}{ccccccc}
   1 &\frac{ k_{1,2}}{k_{1,1}}&\frac{k_{1,3}}{k_{1,1}}& \dots&  \frac{k_{1,\Delta(\overline{G})+1} }{k_{1,1}}& | & \frac{g_1}{k_{1,1}} \label{ms3} \\
    0 & k_{2,2}-\frac{k_{1,2}k_{2,1}}{k_{1,1}}& k_{2,3}-\frac{k_{1,3}k_{2,1}}{k_{1,1}}& \dots&  k_{2,\Delta(\overline{G})+1}-\frac{k_{1,\Delta(\overline{G})+1}k_{2,1}}{k_{1,1}} & |  &-\frac{g_1k_{2,1}}{k_{1,1}} \label {ms5}\\
  \end{array}
\right)
$$
and since  $k_{2,2}\not\in \Q[k_{1,2},k_{2,1},k_{1,1}]$,  $k_{2,2}-\frac{k_{1,2}k_{2,1}}{k_{1,1}}\ne 0$. Hence we  obtain  the echelon  form  of the matrix  dividing  by $\alpha=k_{2,2}-\frac{k_{1,2}k_{2,1}}{k_{1,1}}$
$$
\left(
  \begin{array}{ccccccc}
   1 &\frac{ k_{1,2}}{k_{1,1}}&\frac{k_{1,3}}{k_{1,1}}& \dots&  \frac{k_{1,\Delta(\overline{G})+1} }{k_{1,1}}& | & \frac{g_1}{k_{1,1}} \label{ms3} \\
    0 & 1&\frac{1}{\alpha}( k_{2,3}-\frac{k_{1,3}k_{2,1}}{k_{1,1}})& \dots&\frac{1}{\alpha}(  k_{2,\Delta(\overline{G})+1}-\frac{k_{1,\Delta(\overline{G})+1}k_{2,1}}{k_{1,1}}) & |  &-\frac{g_1k_{2,1}}{\alpha k_{1,1}} \label {ms6}\\
  \end{array}
\right)
$$
Choose $k_{3,3}\ne 0, k_{3,4}\ne 0,\dots,k_{3,\Delta(\overline{G})+1}\ne0$ one  by one  from different  field extensions $\Q[\gamma_3],\dots,\\ \Q[\gamma_{\Delta(\overline{G})+1}]$ such that $\Q[\gamma_i], i=3,\dots,\Delta(\overline{G})+1$ is not   a  field extension in the lattice $L(k_{1,1},\dots,\\ k_{1,\Delta(\overline{G})+1},k_{2,1},\dots,k_{2,\Delta(\overline{G})+1},\gamma_3,\gamma_4,\dots,\gamma_{i-1})$. Therefore, since  $g_1$ is nonzero, we can  choose $g_1\in \Re$ in such a way  that $g_1$  does not  belong  to  any  of  the prior  field extensions used  so far and   satisfy
\begin{eqnarray*}
% \nonumber to remove numbering (before each equation)
 \frac{-g_1 k_{2,1}}{\alpha k_{1,1}}-\frac{1}{\alpha}\left((k_{2,3}-\frac{k_{1,3}k_{2,1}}{k_{1,1}})k_{3,3}+\dots+(k_{2,\Delta(\overline{G})+1}-\frac{k_{1,\Delta(\overline{G})+1}k_{2,1}}{k_{1,1}})k_{3,\Delta(\overline{G})+1}\right) &\ne&0\ \
 \end{eqnarray*}
\noindent which implies
 \begin{eqnarray}
 g_1 &\ne&\frac{-k_{1,1}}{k_{2,1}}\left( (k_{2,3}-\frac{k_{1,3}k_{2,1}}{k_{1,1}})k_{3,2}+\dots+(k_{2,\Delta(\overline{G})+1}-\frac{k_{1,\Delta(\overline{G})+1}k_{2,1}}{k_{1,1}})k_{3,\Delta(\overline{G})+1}\right).\label{cond1}
\end{eqnarray}
\vspace{-1cm}
\noindent Also
\begin{eqnarray*}
\frac{g_1}{k_{1,1}}-\left(\frac{k_{1,2}}{k_{1,1}}k_{3,2}+\frac{k_{1,3}}{k_{1,1}}k_{3,3}+\dots+\frac{k_{1,\Delta(\overline{G})+1}}{k_{1,1}}k_{3,\Delta(\overline{G})+1}\right)&\ne&0
\end{eqnarray*}
\noindent which  implies
\vspace{-1cm}
\begin{eqnarray}
g_1&\ne&k_{1,1}\left(\frac{k_{1,2}}{k_{1,1}}k_{3,2}+\frac{k_{1,3}}{k_{1,1}}k_{3,3}+\dots+\frac{k_{1,\Delta(\overline{G})+1}}{k_{1,1}}k_{3,\Delta(\overline{G})+1}\right).\label{cond2}
\end{eqnarray}
Thus,  choosing $g_1\ne 0$  satisfying \ref{cond1} and \ref{cond2} we get  that  the system \ref{esys1} and \ref{esys2}  is consistent  and at least one  of  its  solutions  satisfies    the adjacency  condition and the orthogonal condition  for  $\overrightarrow{v}_3$ such that none of   the coordinates of  the vectors $\overrightarrow{v_1},\overrightarrow{v_2},\overrightarrow{v_3} \in \Re^{\Delta(\overline{G})+1}$ are   zero.
Therefore   $\{\overrightarrow{v}_1,\overrightarrow{v}_2,\overrightarrow{v}_3\}$ is   an orthogonal  representation in  $\Re^{\Delta(\overline{G})+1}$ for   the  induced  graph $H=K_2\sqcup K_1$ sa\-tisfying  condition (i).
Note  that   if  $\overrightarrow{v}_3=a \overrightarrow{v}_1, a\in \Re$ then  $k_{3,1}= a k_{1,1}$ and $k_{3,\Delta(\overline{G})+1}= a k_{1,\Delta(\overline{G})+1}$  which  implies   that $a =\frac{k_{3,1}}{ k_{1,1}}$ and $a =\frac{k_{3,\Delta(\overline{G})+1}}{ k_{1,\Delta(\overline{G})+1}}$. As a consequence $k_{3,\Delta(\overline{G})+1}=\frac{k_{3,1}}{ k_{1,1}} k_{1,\Delta(\overline{G})+1}$. Hence $k_{3,\Delta(\overline{G})+1}\in \Q[k_{1,1},k_{3,1}, k_{1,\Delta(\overline{G})+1}]$ which  is a contradiction  because  $k_{3,\Delta(\overline{G})+1}$ was  chosen  from a different  field  extension.  Hence,  $\overrightarrow{v}_3$ and  $\overrightarrow{v}_1$ are linearly independent and the vectors $\overrightarrow{v}_1,\overrightarrow{v}_2,\overrightarrow{v}_3$ are pairwise  linearly independent satisfying condition (ii).
\hfill $\Box$
%*****************************************************
%                  MAIN  RESULT
%*****************************************************
\begin{theorem} \label{main} Let $G(V,E)$ be a $\delta$-graph then
$$
\msr(G)\le\Delta(\overline{G})+1=|G|-\delta(G)\label{mrsineq1}
$$
\end{theorem}
\noindent\MakeUppercase{Proof.} Let $G(V,E)$  be a $\delta$ graph. Let $Y_3(\{v_1,v_2,v_3\},E_{Y_3})$ be  the graph induced  by  the  vertices  $v_1,v_2$, and $v_3$ which is  either $3K_1$ or $K_2\sqcup K_1$. By Lemma \ref{lem2}, $Y_3$  has an orthogonal  representation in $\Re^{\Delta(\overline{G})+1}$. Also, from  Lemma \ref{lem2}   we have
\begin{enumerate}
\item [(1)]\label{tcond2} The  components  of vectors  in the orthogonal  representation of  $Y_3$ are all nonzero.
\item [(2)] \label{tcon1} $\overrightarrow{v}\not\in \sp(\overrightarrow{u})$ for each pair   of  distinct  vertices  $u,v$ in  $V_{Y_3}$.
\end{enumerate}
From  the definition  of  $\delta$-graph we have
\begin{enumerate}
\item [(3)] \label{tcond3} $G(V,E)$  can be  constructed starting  with $Y_3(V_{Y_3},E_{Y_3})$ and adding  one  vertex  at a time  such that the newly added  vertex $v_m, m\ge 4$ is adjacent to all prior vertices $v_1,v_2,\dots,v_{m-1}$  except for at most  $\dis{\left\lfloor\frac{m}{2}-1\right\rfloor}$ vertices.
\end{enumerate}
Applying  condition (3) we  get a  sequence  of  subgraphs
$
Y_3,Y_4,\dots,Y_m,\dots,Y_{|G|}\label{seqinG}
$ in $G$  induced by
$
\{v_1,v_2,v_3\},\{v_1,v_2,v_3,v_4\},\dots,\{v_1,v_2,v_3,v_4,\dots,v_m\},$ $\dots, \{v_1,v_2,v_3,v_4,\dots,v_{|G|}\}\label{seqvert}
$
respectively.

We will prove that  $Y_j=(V_j,E_j)$  has an orthogonal  representation  in $\Re^{\Delta(\overline{G})+1}$ for all $j= 3,4,\dots,|G|$, satisfying conditions  (1)  and  (2) above. For  that purpose, consider the orthogonal  representation of  $Y_3$ satisfying  conditions (i) and  (ii)  given  by Lemma \ref{lem2}.
\begin{eqnarray*}
% \nonumber to remove numbering (before each equation)
 \overrightarrow{v_1} &=& k_{1,1}\overrightarrow{e_1}+k_{1,2}\overrightarrow{e_2}+\dots+k_{1,\Delta(\overline{G})+1}\overrightarrow{e}_{\Delta(\overline{G})+1} \label{tm1}\\
 \overrightarrow{v_2} &=& k_{2,1}\overrightarrow{e_1}+k_{2,2}\overrightarrow{e_2}+\dots+k_{2,\Delta(\overline{G})+1}\overrightarrow{e}_{\Delta(\overline{G})+1} \label{tm2} \\
\overrightarrow{v_3} &=& k_{3,1}\overrightarrow{e_1}+k_{3,2}\overrightarrow{e_2}+\dots+k_{3,\Delta(\overline{G})+1} \overrightarrow{e}_{\Delta(\overline{G})+1}\label{tm3}
\end{eqnarray*}
\noindent where  all $k_{i,j}, j=1,2,\dots,\Delta(\overline{G})+1, i=1,2,3$, are nonzero  and are  chosen  from   different field  extensions as in the  proof of Lemma \ref{lem2}. Let $v_4$  be a vertex  of  $G$ such that  $v_4$ is adjacent to all of $v_1,v_2,v_3$ except  at most   $\left \lfloor\frac{4}{2}-1\right\rfloor=1$ vertex. Since  $G$ and  $\overline{G}$  are  simple and  connected   $2\le \Delta(\overline{G})$ and therefore $\Delta(\overline{G}) + 1 \ge 3$.
\begin{claim}\label{claimMainT}
$Y_4$ induced  by $\{v_1,v_2,v_3,v_4\}$ has an orthogonal  representation  in $\Re^{\Delta(\overline{G}) + 1}$  satisfying conditions  (1) and  (2) above.
\end{claim}
\noindent\MakeUppercase{Proof of  claim.}We know  that
$$
d_{\overline{G}}(v_4)\le\Delta(\overline{G})<\Delta(\overline{G})+1
$$
Since $V_{Y_3}=\{v_1,v_2,v_3\}$ from Lemma \ref{lem2} we have  an orthogonal  representation  for the induced  subgraph of $V_{Y_3}$ of  $G$ in $\Re^{\Delta(\overline{G}) + 1}$ satisfying conditions  (1) and  (2). We need  to find a vector $\overrightarrow{v_4}$  for  the vertex $v_4$ where
$$
\overrightarrow{v_4} = k_{4,1}\overrightarrow{e_1}+k_{4,2}\overrightarrow{e_2}+\dots+k_{4,\Delta(\overline{G})+1} \overrightarrow{e}_{\Delta(\overline{G})+1}\label{tm4}
$$
satisfying conditions   (1) and (2).
Since  $v_4$  is adjacent with  all  prior   vertices except  for  at most  one  of  them we  have four  cases:
\begin{enumerate}
\item \label{case1} $v_4$  is adjacent  to $v_1,v_2$ and $v_3$ in $G$.
\item \label{case2} $v_4$  is adjacent  to  only $v_1$ and $v_2$ in $G$.
\item \label{case3} $v_4$  is adjacent  to  only $v_1$ and $v_3$ in $G$.
\item \label{case4} $v_4$  is adjacent  to  only $v_2$ and $v_3$ in $G$.
\end{enumerate}
\noindent {\it Case 1.\ \ \ $v_4\sim v_1,v_4\sim v_2,v_4\sim v_3$ in $G$}.\\
Choose $k_{4,j}, j=1,\dots,\Delta(\overline{G})+1$ as follows:
\begin{enumerate}
\item[] $k_{4,1} =\gamma_{4,1}$ does not  belong to  any of  the  field extensions   in   the lattice  of fields  $L[\Q[k_{i,j}]],\\  i=1,2,3,j=1,2,\dots,\Delta(\overline{G})+1$ .
\item[] $k_{4,2} =\gamma_{4,2}$ does not  belong to  any of  the  field extensions   in   the lattice  of field extensions $L[\Q[k_{i,j},\gamma_{4,1}] ,  i=1,2,3, j=1,2,\dots,\Delta(\overline{G})+1$.
\item[] $k_{4,3} =\gamma_{4,3}$ does not  belong to  any of  the  field extensions   in   the lattice  of field extensions $L[\Q[k_{i,j},\gamma_{4,1},\gamma_{4,2}]] ,  i=1,2,3, j=1,2,\dots,\Delta(\overline{G})+1$.
\item[] Continuing  the process until $k_{4,\Delta(\overline{G})+1}$ to obtain
\item[] $k_{4,\Delta(\overline{G})+1} =\gamma_{4,\Delta(\overline{G})+1}$ does not  belong to  any of  the  field extensions   in   the lattice  of field extensions $L[\Q[k_{i,j},\gamma_{4,1},\gamma_{4,2},\dots,\gamma_{4,\Delta(\overline{G})}]] ,  i=1,2,3, j=1,2,\dots,\Delta(\overline{G})+1$.
\end{enumerate}
Then $k_{4,j} \ne0$ for all $j= 1,2,\dots, \Delta(\overline{G})+1$ and  $\overrightarrow{v_4}\not \in \sp(\overrightarrow{v_i}), i=1,2,3$.
As a consequence $\{\overrightarrow{v_1},\overrightarrow{v_2},\overrightarrow{v_3},\overrightarrow{v_4}\}$ is  an orthogonal  representation of  $Y_4$  at $\Re^{\Delta(\overline{G})+1}$. Note  that if  $\langle\overrightarrow{v}_4,\overrightarrow{v}_1\rangle=0$  then  we  can  solve  this equation  for $k_{4,\Delta(\overline{G})+1}$ which  implies  that  $k_{4,\Delta(\overline{G})+1}\in L[\Q[k_{i,j},\gamma_{4,1},\gamma_{4,2},\dots,\\ \gamma_{4,\Delta(\overline{G})}]] ,  i=1,2,3, j=1,2,\dots,\Delta(\overline{G})+1$ which  is a  contradiction. Therefore, $\langle\overrightarrow{v}_4,\overrightarrow{v}_1\rangle\ne0$.  In  the same  way   we can  prove  that $\langle\overrightarrow{v}_4,\overrightarrow{v}_2\rangle\ne0$  and $\langle\overrightarrow{v}_4,\overrightarrow{v}_3\rangle\ne0$.

\noindent {\it Case 2.\ \ \ $v_4\sim v_1,v_4\sim v_2,v_4\not \sim v_3$ in $G$}.

Since $v_4\sim v_1$ and $v_4\sim v_2$ and $v_4\not \sim v_3$ then
\begin{eqnarray*}
% \nonumber to remove numbering (before each equation)
  \langle\overrightarrow{v_4},\overrightarrow{v_1}\rangle &=& g_{4,1}, \ g_{4,1}\ne 0 \\
  \langle\overrightarrow{v_4},\overrightarrow{v_2}\rangle &=& g_{4,2}, \ g_{4,2}\ne 0\\
  \langle\overrightarrow{v_4},\overrightarrow{v_3}\rangle &=&  0.
\end{eqnarray*}
From these conditions   the    system $S$ in the variables $k_{4,j}, j=1,2,\dots,\Delta(\overline{G})+1$ becomes,
\begin{eqnarray*}
 % \nonumber to remove numbering (before each equation)
  k_{1,1}k_{4,1}+k_{1,2}k_{4,2}+\dots+k_{1,\Delta(\overline{G})+1}k_{4,\Delta(\overline{G})+1} &=&g_{4,1}, g_{4,1} \ne 0\\
   k_{2,1}k_{4,1}+k_{2,2}k_{4,2}+\dots+k_{2,\Delta(\overline{G})+1}k_{4,\Delta(\overline{G})+1} &=& g_{4,2}, g_{4,2} \ne 0 \\
  k_{3,1}k_{4,1}+k_{3,2}k_{4,2}+\dots+k_{3,\Delta(\overline{G})+1}k_{4,\Delta(\overline{G})+1}&=& 0
 \end{eqnarray*}
\noindent where $k_{i,j}, i=1,2,3, j=1,2,\dots,\Delta(\overline{G})+1$ were  chosen   from   different  field extensions as in  the proof  of Lemma \ref{lem2}.
 Since
 $$
 d_{\overline{G}}(v_4)\le\Delta(\overline{G})<\Delta(\overline{G})+1
 $$
 in $\overline{G}$,  the number  of    equations  from  the orthogonal  conditions  in  the systems are at most $\Delta(\overline{G})<\Delta(\overline{G})+1$, which means  that  if  the non-homogeneous system is consistent  then  the  system   will   have  infinitely many  solutions  because   the system  will have at least  one free variable.
 The augmented  matrix  of  the system  becomes
 \begin{center}
 $
 \left(
   \begin{array}{cccccc}
     k_{1,1} &k_{1,2} &\dots & k_{1,\Delta(\overline{G})+1} & | & g_{4,1} \\
    k_{2,1} &k_{2,2} &\dots & k_{2,\Delta(\overline{G})+1} & | & g_{4,2}\\
    k_{3,1} &k_{3,2} &\dots & k_{3,\Delta(\overline{G})+1} & | &0
   \end{array}
 \right)
 $
 \end{center}
 where $g_{4,1}\ne0, g_{4,2}\ne0$  are  nonzero real numbers.

In order to  guarantee  that  all adjacency conditions   are  satisfied consider  the matrix  $3\times(\Delta(\overline{G})+3)$ below in the variables  $k_{4,1},k_{4,2},\dots,k_{4,\Delta(\overline{G})+1},-g_{4,1},-g_{4,2}$. It is possible to  consider $-g_{4,1}$ and $-g_{4,2}$  as   variables because we  only  need  them  to be  nonzero. So  we can  consider them as  two additional  variables of  the homogeneous  system $S_H$ which  has  the  following  augmented matrix:
  \begin{center}
  $
 \left(
   \begin{array}{cccccc}
     k_{1,1} &k_{1,2} &\dots & k_{1,\Delta(\overline{G})+1} & 1 & 0 \\
    k_{2,1} &k_{2,2} &\dots & k_{2,\Delta(\overline{G})+1} & 0 & 1\\
    k_{3,1} &k_{3,2} &\dots & k_{3,\Delta(\overline{G})+1} & 0&0
   \end{array}
 \right).
 $
 \end{center}
Multiplying  the  first row  by  $\frac{1}{k_{1,1}},  k_{1,1}\ne 0$ we get
  \begin{center}
  $
 \left(
   \begin{array}{cccccc}
    1 &\frac{k_{1,2}}{ k_{1,1}} &\dots & \frac{k_{1,\Delta(\overline{G})+1}}{ k_{1,1}} & \frac{1}{ k_{1,1}} & 0 \\
    k_{2,1} &k_{2,2} &\dots & k_{2,\Delta(\overline{G})+1} & 0 & 1\\
    k_{3,1} &k_{3,2} &\dots & k_{3,\Delta(\overline{G})+1} & 0&0
   \end{array}
 \right).
 $
  \end{center}
 Multiplying  the  first  row  by  $-k_{2,1}$ and  adding  the result  to  the second  row and  multiplying  the first  row  by $-k_{3,1}$ and adding  to  the  third   row   we  get
\begin{center}
  $
 \left(
   \begin{array}{cccccc}
    1 &\frac{k_{1,2}}{ k_{1,1}} &\dots & \frac{k_{1,\Delta(\overline{G})+1}}{ k_{1,1}} & \frac{1}{ k_{1,1}} & 0 \\
   0 &k_{2,2}-\frac{k_{1,2}k_{2,1}}{ k_{1,1}} &\dots & k_{2,\Delta(\overline{G})+1}- \frac{k_{1,\Delta(\overline{G})+1}k_{2,1}}{ k_{1,1}} & -\frac{k_{2,1}}{k_{1,1}}& 1\\
   0 &k_{3,2}-\frac{k_{1,2}k_{3,1}}{ k_{1,1}} &\dots & k_{3,\Delta(\overline{G})+1}- \frac{k_{1,\Delta(\overline{G})+1}k_{3,1}}{ k_{1,1}} & 0&0
   \end{array}
 \right).
 $
  \end{center}
Let $\alpha =k_{2,2}-\frac{k_{1,2}k_{2,1}}{ k_{1,1}}$. Since $k_{2,2}\not  \in \Q[k_{1,1},k_{1,2},k_{2,1}]$  by  construction, $\alpha\ne0$  and  we can continue  reducing the matrix   to echelon  form. Then multiplying  the second  row  by $\frac{1}{\alpha}$ we obtain
\begin{center}
  $
 \left(
   \begin{array}{ccccccc}
    1 &\frac{k_{1,2}}{ k_{1,1}}&\frac{k_{1,3}}{ k_{1,1}}  &\dots & \frac{k_{1,\Delta(\overline{G})+1}}{ k_{1,1}} & \frac{1}{ k_{1,1}} & 0 \\
   0 &1 &\frac{1}{\alpha}(k_{2,3}- \frac{k_{1,3}k_{2,1}}{ k_{1,1}})&\dots &\frac{1}{\alpha}( k_{2,\Delta(\overline{G})+1}- \frac{k_{1,\Delta(\overline{G})+1}k_{2,1}}{ k_{1,1}} )& -\frac{k_{2,1}}{\alpha k_{1,1}}&\frac{1}{\alpha}\\
   0 &k_{3,2}-\frac{k_{1,2}k_{3,1}}{ k_{1,1}} & k_{3,3}- \frac{k_{1,3}k_{3,1}}{ k_{1,1}} &\dots & k_{3,\Delta(\overline{G})+1}- \frac{k_{1,\Delta(\overline{G})+1}k_{3,1}}{ k_{1,1}} & 0&0
   \end{array}
 \right).
 $
  \end{center}
  Let $\beta=k_{3,2}-\frac{k_{1,2}k_{3,1}}{ k_{1,1}}$. Multiplying the  second  row  by $-\beta$ and adding  the  result  to  the  third   row  we  get
 \begin{center}
  $
 \left(
   \begin{array}{ccccccc}
    1 &\frac{k_{1,2}}{ k_{1,1}}&\frac{k_{1,3}}{ k_{1,1}}  &\dots & \frac{k_{1,\Delta(\overline{G})+1}}{ k_{1,1}} & \frac{1}{ k_{1,1}} & 0 \\
   0 &1 &\frac{1}{\alpha}(k_{2,3}- \frac{k_{1,3}k_{2,1}}{ k_{1,1}})&\dots &\frac{1}{\alpha}( k_{2,\Delta(\overline{G})+1}- \frac{k_{1,\Delta(\overline{G})+1}k_{2,1}}{ k_{1,1}} )& -\frac{k_{2,1}}{\alpha k_{1,1}}&\frac{1}{\alpha}\\
   0 &0 & \rho_{3,3}&\dots &\rho_{3,\Delta(\overline{G})+1} & \frac{-k_{3,1}}{k_{1,1}}+\frac{k_{2,1}\beta}{\alpha k_{1,1}}&\frac{-\beta}{\alpha}
   \end{array}
 \right)
 $
\end{center}
where  $\rho_{3,3}=k_{3,3}- \frac{k_{1,3}k_{3,1}}{ k_{1,1}}-\frac{\beta}{\alpha}(k_{2,3}- \frac{k_{1,3}k_{2,1}}{ k_{1,1}}) ,\dots,\rho_{3,\Delta(\overline{G})+1}=k_{3,\Delta(\overline{G})+1}- \frac{k_{1,\Delta(\overline{G})+1}k_{3,1}}{ k_{1,1}}-\frac{\beta}{\alpha}( k_{2,\Delta(\overline{G})+1}- \frac{k_{1,\Delta(\overline{G})+1}k_{2,1}}{ k_{1,1}} )$.
Note  that $\rho_{3,3}\ne 0$ otherwise  $k_{3,3}$ belongs  to  a field  extension  of  the lattice
$
L[\Q[k_{1,1},\dots,k_{1,\Delta(\overline{G})+1},k_{2,1},\dots,k_{2,\Delta(\overline{G})+1},k_{3,1},k_{3,2}]]
$
which is a contradiction.

Thus,  multiplying   the  third  row  by $\frac{1}{\rho_{3,3}}$ we  obtain the echelon  form   of  the homogeneous system
$$
\left(
   \begin{array}{cccccccc}
   1 &\frac{k_{1,2}}{ k_{1,1}}&\frac{k_{1,3}}{ k_{1,1}} &\frac{k_{1,4}}{k_{1,1}}&\dots & \frac{k_{1,\Delta(\overline{G})+1}}{ k_{1,1}} & \frac{1}{ k_{1,1}} & 0 \\
0 &1 &\frac{1}{\alpha}(k_{2,3}- \frac{k_{1,3}k_{2,1}}{ k_{1,1}})&\frac{1}{\alpha}(k_{2,4}-\frac{k_{1,4}k_{2,1}}{k_{1,1}})&\dots &\frac{1}{\alpha}( k_{2,\Delta(\overline{G})+1}- \frac{k_{1,\Delta(\overline{G})+1}k_{2,1}}{ k_{1,1}} )& -\frac{k_{2,1}}{\alpha k_{1,1}}&\frac{1}{\alpha}\\
  0 &0 & 1&\frac{\rho_{3,4}}{\rho_{3,3}}&\dots &\frac{\rho_{3,\Delta(\overline{G})+1}}{\rho_{3,3}} &\zeta&\frac{-\beta}{\rho_{3,3}\alpha}
   \end{array}
\right)
$$
where all  the values $\rho_{3,j}\ne0 , j=4,\dots,\Delta(\overline{G})+1$ as well as $\zeta=\frac{-k_{3,1}}{\rho_{3,3}k_{1,1}}+\frac{k_{2,1}\beta}{\rho_{3,3}\alpha k_{1,1}}$ and $\frac{-\beta}{\rho_{3,3}\alpha}$ are  nonzero.
Then   we can choose  values  $k_{4,4},\dots,k_{4,\Delta(\overline{G})+1}, g_{4,1}, g_{4,2}$   nonzero and  chosen   one  by one from diffe\-rent field extensions not in the lattice
$$L[\Q[k_{1,1}],\dots, \Q[k_{1,\Delta(\overline{G})+1}], \Q[k_{2,1}],\dots,\Q[k_{2,\Delta(\overline{G})+1}],\Q[k_{3,1}],\\ \dots, \Q[k_{3,\Delta(\overline{G})+1}]$$
 in  the  following  way.
\begin{enumerate}
\item[] $k_{4,4}\not \in \Q[k_{i,j}],i=1,2,3, j=1,\dots,\Delta(\overline{G})+1$,
\item[] $k_{4,3}\not \in \Q[k_{i,j},k_{4,4}],i=1,2,3, j=1,\dots,\Delta(\overline{G})+1$,
\item[] $\vdots$
\item[] $k_{4,\Delta(\overline{G})+1}\not \in \Q[k_{i,j},k_{4,4},\dots,k_{4,\Delta(\overline{G})}],i=1,2,3, j=1,\dots,\Delta(\overline{G})+1$,
\item[] $g_{4,1}\not \in \Q[k_{i,j},k_{4,4},\dots,k_{4,\Delta(\overline{G})+1}],i=1,2,3, j=1,\dots,\Delta(\overline{G})+1$,
\item[] $g_{4,2}\not \in \Q[k_{i,j},k_{4,4},\dots,k_{4,\Delta(\overline{G})+1}, g_{4,1}],i=1,2,3, j=1,\dots,\Delta(\overline{G})+1$,
\end{enumerate}
Therefore $k_{4,3},k_{4,2}$, and $k_{4,1}$ become
\begin{eqnarray*}
% \nonumber to remove numbering (before each equation)
 k_{4,3} &=&-\frac{\rho_{3,4}k_{4,4}}{\rho_{3,3}}
    +\dots-
   \frac{\rho_{3, \Delta(\overline{G})+1}k_{4, \Delta(\overline{G})+1}}{\rho_{3,3}}+ \frac{k_{3,1}g_{4,1}}{\rho_{3, 3}k_{1,1}}-\frac{k_{2,1}\beta g_{4,1}}{\rho_{3,3}\alpha k_{1,1}} +\frac{\beta g_{4,2}}{\rho_{3,3}\alpha}.\\
 k_{4,2} &=&-(\frac{1}{\alpha}(k_{2,3}- \frac{k_{1,3}k_{2,1}}{ k_{1,1}}))k_{4,3}-\dots-(\frac{1}{\alpha}( k_{2,\Delta(\overline{G})+1}- \frac{k_{1,\Delta(\overline{G})+1}k_{2,1}}{ k_{1,1}} ))k_{4,\Delta(\overline{G})+1}\\
  && +\frac{k_{2,1}g_{4,1}}{\alpha k_{1,1}}-\frac{g_{4,2}}{\alpha}\\
k_{4,1} &=&-\frac{k_{1,2} k_{4,2}}{k_{1,1}}-\dots-
  \frac{k_{1, \Delta(\overline{G})+1}k_{4, \Delta(\overline{G})+1}}{k_{1,1}}-
  \frac{g_{4,1}}{k_{1,1}} \label{k41}
\end{eqnarray*}
Note  that $k_{4,1}$ depends on $k_{4,2}$  and $k_{4,3}$.  Similarly $k_{4,2}$ depends  on  $k_{4,3}$. Therefore  we  should   choose  $k_{4,3}$ first  and  then  back  substitute. But $k_{4,3}$  depends  on  $g_{4,1}$ and  $g_{4,2}$  which  are  free  variables  with  the  restriction  that they  cannot be zero.  Thus,  we  can  choose  $g_{4,1}$  in  a  field   extension  not  in  the lattice  of  $$Q[k_{1,1},\dots,k_{1,\Delta(\overline{G})+1},k_{2,1},\dots,k_{2,\Delta(\overline{G})+1}, k_{3,1},\dots, k_{3,\Delta(\overline{G})+1}, k_{4,4}, \dots,k_{4,\Delta(\overline{G})+1}]
$$
 and   $g_{4,2}$   in  some  field  extension not in   the lattice  of $$Q[k_{1,1},\dots,k_{1,\Delta(\overline{G})+1},k_{2,1},\dots,k_{2,\Delta(\overline{G})+1}, k_{3,1},\dots, k_{3,\Delta(\overline{G})+1},\\ k_{4,4},\dots, k_{4,\Delta(\overline{G})+1},g_{4,1}]$$
 in  such  a  way   that  $k_{4,3},k_{4,2}$  and  $k_{4,1}$ are  nonzero.

As a consequence,  the  system $S$ is consistent  and  there exist a solution  of   values $k_{4, j}\ne 0 , \  j=1,2,\dots, \Delta(\overline{G})+1$    chosen  from different  field  extensions  such that  the vector $\overrightarrow{v_4}$ satisfies  all  adjacency conditions and  orthogonal conditions  with  the vectors $\overrightarrow{v_1},\overrightarrow{v_2},\overrightarrow{v_3}$. Then  the vector  $\overrightarrow{v_4}=\sum_{j=1}^{\Delta(\overline{G})+1} k_{4,j} \overrightarrow{e_j}$ satisfies  conditions  (1) and  (2) and  therefore $\{\overrightarrow{v_1},\overrightarrow{v_2},\overrightarrow{v_3},\overrightarrow{v_4}\}$ is an orthogonal  representation  of $Y_4$ in $\Re^{\Delta(\overline{G})+1}$.\\
\noindent {\it Case 3 .\ \ \ $v_4\sim v_1,v_4\not \sim v_2,v_4 \sim v_3$ in $G$}.\\
From  the adjacency  conditions $v_4\sim v_1,v_4 \sim v_3$  and  orthogonal condition $v_4\not \sim v_2$ we  get the equations:
\begin{eqnarray*}
% \nonumber to remove numbering (before each equation)
  \langle\overrightarrow{v_1},\overrightarrow{v_4}\rangle &=& g_{4,1}, g_{4,1}\ne 0 \\
  \langle\overrightarrow{v_2},\overrightarrow{v_4}\rangle &=& 0\\
  \langle\overrightarrow{v_3},\overrightarrow{v_4}\rangle&=& g_{4,2}, g_{4,2}\ne 0
\end{eqnarray*}
Interchanging  the  second  and   third equations  we get a  system $S$ similar to  case  2. Since  all   the scalars $k_{i,j}, i= 1,2,3 ,  j=1,2,\dots,\Delta(\overline{G})+1 $    are not  zero and  were  chosen    from  different  field  extensions  the  same  reasoning as in  case 2 applies  and     the conclusion  holds   for case  3.\\
\noindent {\it Case 4.\ \ \ $v_4\not \sim v_1,v_4\sim v_2,v_4 \sim v_3$ in $G$}.

From  the adjacency  conditions $v_4\sim v_2,v_4 \sim v_3$  and  orthogonal condition $v_4\not \sim v_1$ we  get  the equations:
\begin{eqnarray*}
% \nonumber to remove numbering (before each equation)
  \langle\overrightarrow{v_1},\overrightarrow{v_4}\rangle &=& 0 \\
 \langle\overrightarrow{v_2},\overrightarrow{v_4}\rangle &=& g_{4,2}, g_{4,2}\ne 0\\
  \langle\overrightarrow{v_3},\overrightarrow{v_4}\rangle&=&g_{4,3}, g_{4,3}\ne 0.
\end{eqnarray*}
Interchanging  the  first  and  the  third equations  we get a  system $S$ similar to  case  2. Since  all   the scalars $k_{i,j}, i= 1,2,3 ,  j=1,2,\dots,\Delta(\overline{G})+1 $    are not  zero and  were  chosen    from  different  field  extensions  the  same  reasoning as in  case 2 applies  and the conclusion  holds   for case  4.

As a consequence, in all  of the cases   we  get an orthogonal  representation  for $Y_4$  in $\Re^{\Delta(\overline{G})+1}$ satisfying   the conditions  (1) and  (2). This completes  of the proof of the claim \ref{claimMainT}.

Assume that  for  any $Y_{m-1}=(V_{Y_{m-1}},E_{Y_{m-1}}), V_{Y_{m-1}}= \{v_1,v_2,\dots,v_{m-1}\}$
it is possible  to get  an orthogonal representation of  $Y_{m-1}$ in
$\Re^{\Delta(\overline{G})+1}$. Let $\overrightarrow{v}_1,\overrightarrow{v}_2,\dots,\overrightarrow{v}_{m-1}$ be of  the form

 $$
 \overrightarrow{v}_i=\dis{\sum_{j=1}^{\Delta(\overline{G})+1}k_{i,j}\overrightarrow{e}_j}
 $$
satisfying  conditions  (1) and (2)  where $k_{i,j}\ne0 $ for all $i=1,2,\dots,m , j=1,2,\dots,\Delta(\overline{G})+1$ , chosen  from different  field  extensions.

We need  to  prove   that if $v_{m}$  is adjoined to $Y_{m-1}$  to get $Y_m$ such that $v_{m}$ is adjacent to all prior  vertices  except at most $\left\lfloor\frac{m}{2}-1\right\rfloor$ vertices then  $Y_{m}$ has an orthogonal  representation of  vectors
$ \overrightarrow{v}_1,\overrightarrow{v}_2,\dots,\overrightarrow{v}_{m}$ in $\Re^{\Delta(G)+1}$ satisfying  conditions (1) and  (2).
Assume  that  $v_{m}$  has  an associated  vector $\overrightarrow{v}_m$ such  that
$$
\overrightarrow{v}_{m}= k_1^{m}\overrightarrow{e_1}+ k_2^{m}\overrightarrow{e_2}+\dots+ k_{\Delta(\overline{G})+1}^{m}\overrightarrow{e}_{\Delta(\overline{G})+1}.
$$
 The  vertex $v_{m}$ is adjacent  to all  prior  vertices  $ v_1,v_2,\dots,v_{m-1}$  except at most  $t\le\left\lfloor\frac{m}{2}-1\right\rfloor$ vertices in $G$. Then   we see   that   $\overrightarrow{v}_{m}$   satisfies at least $m-1-t$ adjacency  conditions    and   $t$ orthogonal conditions.

 Let $\rho$ be a permutation of $(1,2,\dots,m-1)$. Suppose $v_{\rho(1)},v_{\rho(2)},\dots,v_{\rho(m-1-t)}$ are  adjacent to $v_m$ and $v_{\rho(m-t)},
 v_{\rho(m-t+1)}, \dots,v_{\rho(m-2)},v_{\rho(m-1)}$  are  not  adjacent   to  $v_{m}$. The vectors  $\overrightarrow{v}_{\rho(1)},\overrightarrow{v}_{\rho(2)},\dots,\\
 \overrightarrow{v}_{\rho(m-1-t)}, \overrightarrow{v}_{\rho(m-t)},\overrightarrow{v}_{\rho(m-t+1)},\dots,\overrightarrow{v}_{\rho(m-1)}$ and $\overrightarrow{v}_{m}$  satisfy  the  system  $S$ given by:
\begin{eqnarray*}
% \nonumber to remove numbering (before each equation)
 \langle\overrightarrow{v}_{\rho(1)},\overrightarrow{v}_{m}\rangle &=&g_{m,1},\ \  g_{m,1}\ne 0\\
  \langle\overrightarrow{v}_{\rho(2)},\overrightarrow{v}_{m}\rangle &=&g_{m,2},\ \  g_{m,2}\ne 0\\
   \vdots&\vdots&\vdots \\
   \langle\overrightarrow{v}_{\rho(m-1-t)},\overrightarrow{v}_{m}\rangle &=&g_{m,m-1-t},\ \  g_{m,m-1-t}\ne 0\\
  \langle\overrightarrow{v}_{\rho(m-t)},\overrightarrow{v}_{m}\rangle&=& 0\\
  \langle\overrightarrow{v}_{\rho(m-t+1)},\overrightarrow{v}_{m}\rangle&=& 0 \\
\vdots &\vdots& \vdots\\
\langle\overrightarrow{v}_{\rho(m-1)},\overrightarrow{v}_{m}\rangle &=&0
\end{eqnarray*}
\noindent containing $m-1-t$  equations from  the adjacency conditions  and $t$  equations from  the  orthogonal conditions.
Since  the vector $\overrightarrow{v}_{\rho(i)}, i=1,2,\dots,m-1$ has  the form
\begin{equation}\label{vectorm}
   \overrightarrow{v}_{\rho(i)}= k_{\rho(i),1}\overrightarrow{e}_1+k_{\rho(i),2}\overrightarrow{e}_2+\dots+k_{\rho(i),\Delta(\overline{G})+1}\overrightarrow{e}_{\Delta(\overline{G})+1}
\end{equation}
where all  $k_{\rho(i),j}, i=1,2,\dots,m-1 ,  j= 1,2,\dots,\Delta(\overline{G})+1$   are not  zero  and chosen  from different  field extensions, the system $S$ has  the  form:
\begin{eqnarray*}
% \nonumber to remove numbering (before each equation)
  k_{\rho(1),1}k_{m,1}+ k_{\rho(1),2}k_{m,2}+\dots+ k_{\rho(1),\Delta(\overline{G})+1}k_{m,\Delta(\overline{G})+1} &=& g_{m,1} \\
 k_{\rho(2),1}k_{m,1}+ k_{\rho(2),2}k_{m,2}+\dots+ k_{\rho(2),\Delta(\overline{G})+1}k_{m,\Delta(\overline{G})+1} &=&g_{m,2} \\
 \vdots &\vdots& \vdots \\
  k_{\rho(m-1-t),1}k_{m,1}+ k_{\rho(m-1-t),2}k_{m,2}+\dots+ k_{\rho(m-1-t),\Delta(\overline{G})+1}k_{m,\Delta(\overline{G})+1} &=&g_{m,m-1-t}   \\
  k_{\rho(m-t),1}k_{m,1}+ k_{\rho(m-t),2}k_{m,2}+\dots+ k_{\rho(m-t),\Delta(\overline{G})+1}k_{m,\Delta(\overline{G})+1}&=& 0\\
  k_{\rho(m-t+1),1}k_{m,1}+ k_{\rho(m-t+1),2}k_{m,2}+\dots+ k_{\rho(m-t+1),\Delta(\overline{G})+1}k_{m,\Delta(\overline{G})+1} &=& 0\\
 \vdots &\vdots& \vdots \\
 k_{\rho(m-1),1}k_{m,1}+ k_{\rho(m-1),2}k_{m,2}+\dots+ k_{\rho(m-1),\Delta(\overline{G})+1}k_{m,\Delta(\overline{G})+1}&=&0
\end{eqnarray*}
where  $g_{m,1}\ne0,g_{m,2}\ne0, \dots,g_{m,m-1-t}\ne0$.
Since  $g_{m,1},g_{m,2}, \dots,g_{m,m-1-t}$  could be  any  nonzero real numbers   satisfying  the adjacency conditions we can consider   them  as additional  $m-1-t$ variables  under  the restriction   that they  cannot  be  zero. Therefore we can consider  a    homogeneous  system  $S_H$ of  $m-1$  equations   in  $m-t+\Delta(\overline{G})$ variables $k_1^{m},k_2^{m}, \dots,k_{\Delta(\overline{G})+1}^{m}, -g_{m,1},-g_{m,2},\\ \dots,-g_{m,m-1-t}$. Now, by hypothesis $t\le \left\lfloor\frac{m}{2}-1\right\rfloor$. Since  $t\le d_{\overline{G}}(v_{m})\le \Delta(\overline{G})<\Delta(\overline{G})+1$ the homogeneous system $S_H$  contains   at least  one more variable   than   the number  of  equations. Hence the  system $S_H$ given by

\begin{eqnarray*}
k_{\rho(1),1}k_{m,1}+ k_{\rho(1),2}k_{m,2}+\dots+ k_{\rho(1),\Delta(\overline{G})+1}k_{m,\Delta(\overline{G})+1}+(-g_{m,1}) &=&0 \\
 k_{\rho(2),1}k_{m,1}+ k_{\rho(2),2}k_{m,2}+\dots+ k_{\rho(2),\Delta(\overline{G})+1}k_{m,\Delta(\overline{G})+1}+(-g_{m,2})&=&0 \\
 \vdots &\vdots& \vdots\\
 k_{\rho(m-1-t),1}k_{m,1}+ k_{\rho(m-1-t),2}k_{m,2}+\dots+ k_{\rho(m-1-t),\Delta(\overline{G})+1}k_{m,\Delta(\overline{G})+1}+(- g_{m,m-1-t})&=&0   \\
  k_{\rho(m-t+),1}k_{m,1}+ k_{\rho(m-t),2}k_{m,2}+\dots+ k_{\rho(m-t),\Delta(\overline{G})+1}k_{m,\Delta(\overline{G})+1}&=& 0 \\
  k_{\rho(m-t+1),1}k_{m,1}+ k_{\rho(m-t+1),2}k_{m,2}+\dots+ k_{\rho(m-t+1),\Delta(\overline{G})+1}k_{m,\Delta(\overline{G})+1} &=& 0\\
 \vdots &\vdots& \vdots \\
 k_{\rho(m-1),1}k_{m,1}+ k_{\rho(m-1),2}k_{m,2}+\dots+ k_{\rho(m-1),\Delta(\overline{G})+1}k_{m,\Delta(\overline{G})+1}&=&0
\end{eqnarray*}

\noindent has infinitely many  solutions. Therefore,  it is enough   to  show   that there exist at least one  solution   for  $S_H$ satisfying  the condition that none of the $k_{m,1},k_{m,2},\dots,k_{m,\Delta(\overline{G})+1},g_{m,1},\\ \dots,g_{m,m-1-t}$ are  zero. This implies that  the system $S$ has a solution  which  satisfies all adjacency conditions, all orthogonal conditions, and  conditions (1) and (2). For that purpose consider    the $(m-1)\times (m-t+\Delta(\overline{G}))$ matrix $A$  of  the homogeneous system given  on  the next  page.

Let $\overrightarrow{g}=(-g_{m,1},-g_{m,2},\dots,-g_{m,m-1-t})^T$. We consider  the  two  cases  where $m-1\le \Delta(\overline{G})+1$ and  $m-1> \Delta(\overline{G})+1 $.

\noindent {\it  Case 1.}\ \   $m-1\le \Delta(\overline{G})+1$

In this case  the number  of equations in the non-homogeneous system $S$ is at most  the number  of unknowns $k_{m,1},k_{m,2},\dots,k_{m,\Delta(\overline{G})+1}$.
\vskip -0.5cm
$$
A=\left(\begin{array}{ccccccccccc}
   k_{\rho(1),1}& k_{\rho(1),2} & k_{\rho(1),3} &\dots & k_{\rho(1),\Delta(\overline{G})+1} &|& 1 & 0 & \dots& 0 &0 \\
 k_{\rho(2),1} &k_{\rho(2),2} &k_{\rho(2),3} & \dots & k_{\rho(2),\Delta(\overline{G})+1} &|&0 &1 & \dots&0& 0 \\
 k_{\rho(3),1} & k_{\rho(3),2} & k_{\rho(3),3} & \dots &k_{\rho(3),\Delta(\overline{G})+1} &|&0 & 0 & 1 &\vdots &\vdots \\
\vdots& \vdots & \vdots & \vdots &\vdots & \vdots & \vdots &\vdots&\vdots &\vdots& \vdots \\
 k_{\rho(m-t-2),1} & k_{\rho(m-t-2),2} & k_{\rho(m-t-2),3} & \dots& k_{\rho(m-t-2),\Delta(\overline{G})+1} &|& 0 & 0 &0 & 1 &0 \\
 k_{\rho(m-t-1),1} & k_{\rho(m-t-1),2}& k_{\rho(m-t-1),3} & \dots& k_{\rho(m-t-1),\Delta(\overline{G})+1}&|&0 & 0 &0 &0 & 1 \\
k_{\rho(m-t),1} & k_{\rho(m-t),2}& k_{\rho(m-t),3} & \dots& k_{\rho(m-t),\Delta(\overline{G})+1} &|&0 & 0 &0 &0 & 0 \\
 \vdots & \vdots &\vdots& \vdots & \vdots & \vdots & \vdots & \vdots & \vdots &\vdots& \vdots \\
k_{\rho(m-1),1} & k_{\rho(m-1),2} &k_{\rho(m-1),3} & \dots & k_{\rho(m-1),\Delta(\overline{G})+1}&| &0 & 0 & 0 & 0 & 0
\end{array}\right)
$$
Thus  $A$ can be  row  reduced  to one of  the  following  two  echelon form written in  block form:
\begin{enumerate}
\item[I.]$$
B= \left(B_1|B_2\right)={\tiny \left(\begin{array}{cccccccccccccccccc}
  1&\ast & \ast &\dots &\ast&\ast&\ast&\dots&\ast&\ast&\dots&\ast&| & \delta_1 & 0 & \dots& 0 &0 \\
0 &1 &\ast & \dots &\ast&\ast&\ast&\dots&\ast&\ast&\dots& \ast&| &\ast &\delta_2 & \dots&0& 0 \\
0 & 0 & 1 & \dots&\ast&\ast &\ast&\dots&\ast&\ast&\dots& \ast&| &\ast&\ast& \delta_3 &0&0 \\
\vdots& \vdots & \vdots &\dots& \vdots &\vdots& \vdots &\dots&\vdots&\vdots&\dots& \vdots&|&\vdots&\vdots&\vdots &\vdots& \vdots \\
0 & 0 &0 & \dots&1&\ast&\ast&\dots&\ast&\ast&\dots&\ast&| &\ast &\ast &\dots& \delta_{m-2-t} &0 \\
0 &0& 0&\dots &0&1&\ast&\dots&\ast&\ast&\dots&\ast&| &\ast & \ast &\dots &\ast& \delta_{m-1-t}  \\
0 &0& 0&\dots &0&0&1&\dots&\ast&\ast&\dots&\ast&| &\ast & \ast &\dots &\ast& \ast  \\
\vdots& \vdots & \vdots &\dots& \vdots &\vdots& \vdots &\dots&\vdots&\vdots&\dots& \vdots&|&\vdots&\vdots&\vdots &\vdots& \vdots \\
0 &0& 0&\dots &0&0&0&\dots&1&\ast&\dots&\ast&| &\ast & \ast &\dots &\ast& \ast  \\
\end{array}\right)}
$$
\item[II.] $$
B=\left(B_3|B_2\right)={\tiny \left(\begin{array}{cccccccccccc}
  1&\ast & \ast &\dots &\ast&\ast&| & \delta_1 & 0 & \dots& 0 &0 \\
0 &1 &\ast & \dots &\ast&\ast&| &\ast &\delta_2 & \dots&0& 0 \\
0 & 0 & 1 & \dots&\ast&\ast &| &\ast&\ast& \delta_3 &0&0 \\
\vdots& \vdots & \vdots &\dots& \vdots &\vdots&|&\vdots&\vdots&\vdots &\vdots& \vdots \\
0 & 0 &0 & \dots&1&\ast&| &\ast &\ast &\dots& \delta_{m-t-1} &0 \\
0 &0& 0&\dots &0&1&| &\ast & \ast &\dots &\ast& \delta_{m-t}  \\
\end{array}\right)}.
$$
\end{enumerate}
The  matrix $B_1$ is a block of  size $(m-1)\times (\Delta(\overline{G})+1)$ where  $m-1< \Delta(\overline{G})+1$.  The matrix $B_2$ is a block matrix of  size $(m-1)\times (m-1-t)$. Matrix  $B_3$  is  a square matrix of  size  $m-1 (=\Delta(\overline{G})+1)$. In  these  blocks $\ast$ denotes a  nonzero entry.

Suppose  matrix  $B$   is  of   type I. For  each  vector $\overrightarrow{v}_i, i=4,5,\dots,m-1$  the entries are found  in  field  extensions  which  are  not  in  the lattice  of  the previous  field extensions.

In  the  block  matrix $B_2$  we  have:
\begin{itemize}
\item [] $[B_2]_{1,1}=\delta_1= \frac{1}{k_{\rho(1),1}}\ne 0$,
\item[] $[B_2]_{2,2}=\delta_2= \left(k_{\rho(2),2}-\frac{k_{\rho(1),2}\cdot k_{\rho(2),1}}{k_{\rho(1),1}}\right)^{-1}\ne 0$.
\item []$[B_2]_{3,3}=\delta_3=\frac{1}{\alpha}\ne 0,  \alpha\in Q[k_{\rho(1),1},k_{\rho(1),2},k_{\rho(1),3},k_{\rho(2),1},k_{\rho(2),2},k_{\rho(2),3},k_{\rho(3),1},k_{\rho(3),2},k_{\rho(3),3}]$.
\end{itemize}
 Continuing  this process   we  get  that   all  the entries  on  the  diagonal  of  $B_2$ to be nonzero.  So  all  the  rows  of $B_2$  have  at least  one  entry nonzero.Thus,
 \begin{itemize}
 \item[] $\left[B_2\overrightarrow{g}\right]_1=-\delta_1\cdot g_{m,1}= \frac{-g_{m,1}}{k_{1,1}}$.  Choosing  $g_{m,1}$  not  in the lattice   generated  by  the  previous  field  extensions  for  $k_{i,j}, i=1,2,\dots,m-1 ,  j=1,2,\dots,\Delta(\overline{G})+1$,  we get  $\left[B_2\overrightarrow{g}\right]_1\ne 0$.
 \item[] $\left[B_2\overrightarrow{g}\right]_2=   \alpha_{2,1}g_{m,1}+\delta_2 g_{m,2},  \alpha_{2,1}\ne 0, \delta_2\ne 0$.  Then    we can  choose  a value for  $g_{m,2}$  from   a  field extension  not   in  the lattice of  fields generated  by   the previous values $k_{i,j}, i=1,2,\dots,m-1, j=1,2,\dots, \Delta(G)+1,g_{m,1}$ so that  $\left[B_2\overrightarrow{g}\right]_2\ne 0$.
 \item[]  $\left[B_2\overrightarrow{g}\right]_3=   \alpha_{3,1}g_{m,1}+\alpha_{3,2}g_{m,2}+\delta_3 g_{m,3},  \alpha_{3,1}\ne 0, \delta_3\ne 0$. As above we can   choose  $g_{m,3}\ne 0$ and  such  that $\left[B_2\overrightarrow{g}\right]_3\ne 0$ by  taking  $g_{m,3}$ neither in $Q[\frac{\alpha_{3,1}g_{m,1}+\alpha_{3,2}g_{m,2}}{-\delta_3}]$ nor  in any  of  the  previous  field  extensions.
 \end{itemize}
Continuing  this   process and applying similar choices  we see  that $g_{m,4},g_{m,5}\dots,g_{m,m-1}$ can  be  chosen nonzero.
Moreover,  matrix $B_1$ shows  that  there is at least  one  free variable for  the solution  of $k_{m,j}, j=1,2,\dots,\Delta(\overline{G})+1$.  All of  these  free  variables  can be chosen  in different  field extensions  such that   all   other  unknowns are not  zero. Otherwise  it is possible  to show   that  the last choice  belongs  to  the  field  containing  all  the  previous chosen  values   which  is a contradiction.  Now,  suppose  that $k_{m,1},k_{m,2},\\ \dots, k_{m,r},  r<\Delta(\overline{G})+1$  can be written   in terms  of  $k_{m,r+1},\dots,k_{m,\Delta(\overline{G}+1)}$ as
\begin{equation}\label{formofks}
k_{m,i}= \alpha_{r+1}k_{m,r+1}+\dots+\alpha_{\Delta(\overline{G})+1}k_{m,\Delta(\overline{G})+1}+ \varphi_i(g_{m,1},g_{m,2},\dots, g_{m,m-1}) , i=1,2,\dots,r
\end{equation}
where $ \alpha_{r+1},\dots,\alpha_{\Delta(\overline{G})+1}$ are all nonzero and $\varphi_i$ is a linear  combination  of  $g_{m,1},\dots,g_{m,m-1}$. Also $\varphi_i(g_{m,1},g_{m,2},\dots, g_{m,m-1})\ne 0 ,i=1,2,\dots,r $. Thus,by  choosing   values\\ $k_{m,r+1},\dots,k_{m,\Delta(\overline{G})+1}$ in  different  field  extensions    and substituting  them  in  \ref{formofks},  we  obtain   that $k_{m,i}\ne 0, i=1,2,\dots,r$.

As a consequence, the  vector $\overrightarrow{v}_m$ exists  and all of  its  entries are nonzero.

If matrix $B$  is of type II we apply  same process as in  case of type I. Again,  we can  obtain  the  vector  $v_m$ having  all its  entries  nonzero and $g_{m,i}\ne 0$  for $i=1,\dots,m-1$.

\noindent {\it  Case 2.}\ \ $m-1 > \Delta(\overline{G})+1$

In  this  case  the number of equations in the non-homogeneous system $S$ is more  than  the number of  unknowns $k_{m,1},k_{m,2},\dots,k_{m,\Delta(\overline{G})+1}$.

We need  to  analyze   three possible  subcases where  $m-1-t< \Delta(\overline{G})+1,  m-1-t=\Delta(\overline{G})+1$ or $m-1-t>\delta(\overline{G})+1$.
\begin{enumerate}
\item  If  $m-1-t< \Delta(\overline{G})+1$ then  matrix $B$   has the  form
$$
B= {\tiny \left(\begin{array}{cccccccccccccccccc}
  1&\ast & \ast &\dots &\ast&|&\ast&\dots&\ast&| &\ast &0&0&0& 0 & \dots& 0 &0 \\
0 &1 &\ast & \dots &\ast&|&\ast&\dots&\ast&| &\ast &\ast& 0&0& 0& \dots&0& 0 \\
0 & 0 & 1 & \dots&\ast&|&\ast&\dots&\ast&| &\ast&\ast& \ast &0&0&\dots&0&0 \\

\vdots& \vdots & \vdots&\dots& \vdots &|&\vdots&\dots&\vdots&|& \vdots &\dots& \vdots&\vdots&\vdots&\vdots &\vdots& \vdots \\
 0 & 0 &0 & \dots&1&|&\ast&\dots&\ast&|&\ast&\dots&\ast&\ast&\ast &\dots& \ast &\ast \\
  -&-&-&-&-&-&-&\dots&-&|&-&-&-&-&-&\dots&-&- \\
  0 & 0 &0 & \dots&0&|&1&\dots&\ast&|&\ast&\dots&\ast &\ast &\ast &\dots& \ast &\ast \\
 \vdots & \vdots &\vdots&\vdots&\vdots&|& \vdots&\dots & \vdots &|&\vdots& \dots&\vdots & \vdots & \vdots & \vdots &\vdots& \vdots \\
0 & 0 &0 & \dots&0&|&0&\dots&\ast&|&\ast&\dots&\ast &\ast &\ast &\dots& \ast &\ast \\
0 &0& 0&\dots &0&|&0&\dots&1&|&\ast&\dots&\ast &\ast & \ast &\dots &\ast& \ast  \\
 %-&-&-&-&-&-&-&\dots&-&|&-&-&-&-&-&\dots&-&- \\
0 & 0&0 & \dots&0&|&0&\dots& 0&|&\ast&\dots&\ast &\ast & \ast &\dots &\ast & \ast \\
 \vdots & \vdots &\vdots&\vdots&\vdots&|& \vdots&\dots & \vdots &|&\vdots& \dots&\vdots & \vdots & \vdots & \vdots &\vdots& \vdots \\
0 & 0 &0 & \dots &0&|&0&\dots&0&|&\ast&\dots&\ast&\ast &\ast & \dots& \ast & \ast
\end{array}\right)}
$$
$$
B=\left(
  \begin{array}{ccc}
    B_1 & B_2 & B_3 \\
    B_4& B_5& B_6 \\
  \end{array}
\right)
$$

The  matrix $B_1$ is a square matrix  $(m-1-t)\times (m-1-t)$, matrix $B_2$ has  size $(m-1-t)\times (\Delta(\overline{G})+ 2 + t-m)$, matrix  $B_3$ is a square matrix  of  size  $m-1-t$.  Matrix $B_4$ is a zero matrix  of  size $t\times m-1-t$.
The  central  blocks $B_2, B_5$   form  a  block  of  size $(m-1)\times (\Delta(\overline{G})+2+ t-m)$ and  corresponds   to the columns  of   free  variables $k_{m,m-t}, \dots, k_{m,\Delta(\overline{G})+1}$ of  the system $S$. The  block  $B_6$ has size $t\times (m-1-t)$ .

  Consider  the block matrix $\left(B_5\ B_6\right)$ of  size $t\times(\Delta(\overline{G})+1)$. Recalling  that the   value $t$  is  the  number of  orthogonal conditions  for  $v_m$ in $G$ which is  equivalent  to  $\d_{\overline{G}}(v_m)$ we  get  $t\le\Delta(\overline{G})<\Delta(\overline{G})+1$. As a consequence,  the   homogeneous system $\left(B_5\ B_6\right)\overrightarrow{w}=0$  where $\overrightarrow{w}$ is  a vector  of  size $(\Delta(\overline{G})+1)\times 1$  in  the variables $k_{m,m-t},\dots, k_{m,\Delta(\overline{G})+1},(-g_{m,1}),\dots,\\ (-g_{m,m-1-t})$,  has  infinitely many  solutions  depending on at least one  free  variable. Choosing  these  free  variables  in  different  field  extensions  as  we  did  previously,  we  get  nonzero values   for   $k_{m,m-t},\\ \dots,k_{m,\Delta(\overline{G})+1},g_{m,1},\dots,g_{m,m-1-t}$. We  get  the  values  of  the remaining  unknowns $k_{m,i}, i=1,2,\dots,m-1-t$ of  the system $S$ applying  back  substitution. Since  all  the entries   with $\ast$ in  the  block $B_1$ are nonzero  and  belong  to different  field  extensions, the  values $k_{m,i}, i=1,2,\dots,m-1-t$ are also  nonzero.

As  a consequence,  the  non-homogeneous system  $S$ is  consistent  and  the vector $v_m$ with no  zero entries exists.

\item  If $m-1-t= \Delta(\overline{G})+1$ then  the matrix  $B$  has  the form
\begin{eqnarray*}
B&=& \left(
    \begin{array}{ccc}
      B_1&| & B_2 \\
      -&&-\\
      0 &|& R \\
    \end{array}
  \right)\\
  &=&{\tiny \left(\begin{array}{cccccccccccccccc}
  1&\ast & \ast &\dots &\ast&\dots&\ast&| &\delta_1 &0&0&0& 0 & \dots& 0 &0 \\
0 &1 &\ast & \dots &\ast&\dots&\ast&| &\ast &\delta_2& 0&0& 0& \dots&0& 0 \\
0 & 0 & 1 & \dots&\ast&\dots&\ast&| &\ast&\ast&\delta_3 &0&0&\dots&0&0 \\
\vdots& \vdots & \vdots&\dots& \vdots &\dots&\vdots&|& \vdots &\dots& \vdots&\vdots&\vdots&\vdots &\vdots& \vdots \\
 0 & 0 &0 & \dots&1&\dots&\ast&|&\ast&\dots&\ast&\ast&\ast&\dots& 0 &0 \\
0 & 0 &0 & \dots&0&\dots&\ast&|&\ast&\dots&\ast &\ast &\ast &\dots& \delta_{m-t-2} &0 \\
0 &0& 0&\dots &0&\dots&1&|&\ast&\dots&\ast &\ast & \ast &\dots &\ast& \delta_{m-1-t}  \\
 -&-&-&-&-&\dots&-&|&-&-&-&-&-&\dots&-&- \\
0 & 0&0 & \dots&0&\dots& 0&|&1&\dots&\ast &\ast & \ast &\dots &\ast & \ast \\
 \vdots & \vdots &\vdots&\vdots& \vdots&\dots & \vdots &\vdots&\vdots& \dots&\vdots & \vdots & \vdots & \vdots &\vdots& \vdots \\
0 & 0 &0 & \dots &0&\dots&0&|&0&\dots&1&\ast &\ast & \dots& \ast & \ast
\end{array}\right)}.
\end{eqnarray*}
In  this  case the matrices  $B_1$ and $B_2$  are  square matrices  of  size $m-1-t (=\Delta(\overline{G})+1)$. The  matrix $R$  has  size $t\times m-1-t$.  Since  $t < \left \lfloor\frac{m}{2}-1\right \rfloor$   we  get   that  $2t<m-2< m-1$ which  implies that  $t<m-1-t$.

Therefore   the system $R\overrightarrow{g}=0$ has infinitely many  solutions  with $m-1-2t$ free variables.  Taking   the free variables   from  different  field  extensions we get  all the values  $g_{m,1}\dots,g_{m,m-1-t}$  nonzero. Substituting $g_{m,1}\dots,g_{m,m-1-t}$  in  the equations   of  the  system  $B_2\overrightarrow{g}=0$  we get $[B_2\overrightarrow{g}]_i \ne 0$ for all $i=1,\dots,m-1-t$. Otherwise, it is possible  to  show  that  the last  choice belongs  to  the  field  containing  all previous  chosen values    which  is a contradiction. This implies   that $k_{m,\Delta(\overline{G})+1}=[B_2\overrightarrow{g}]_{m-1-t}\ne 0$ from  the last  row of $(B_1\ B_2)$.

Applying   back  substitution  and  similar argument  with  $g_{m,1}\dots,g_{m,m-1-t}$ we conclude  that $k_{m,1}, \dots,k_{m,\Delta(\overline{G})}$ are also  nonzero.
 Thus the system $S$ has a solution  with nonzero  values for  the unknowns. As a consequence there  exists a  vector $\overrightarrow{v}_m$  satisfying  all  the adjacency  conditions and orthogonal conditions.

\item $m-1-t> \Delta(\overline{G})+1$  then  matrix  $B$ has  the form
%\vspace*{-2cm}
$$
B= \left(
    \begin{array}{ccc}
      B_1&| & B_2 \\
      -&&-\\
      0 &|& R \\
    \end{array}
  \right)={\tiny \left(\begin{array}{cccccccccccccccc}
  1&\ast & \ast &\dots &\ast&\dots&\ast&| &\ast &0&0&0& 0 & \dots& 0 &0 \\
0 &1 &\ast & \dots &\ast&\dots&\ast&| &\ast &\ast& 0&0& 0& \dots&0& 0 \\
0 & 0 & 1 & \dots&\ast&\dots&\ast&| &\ast&\ast& \ast &0&0&\dots&0&0 \\

\vdots& \vdots & \vdots&\dots& \vdots &\dots&\vdots&|& \vdots &\dots& \vdots&\vdots&\vdots&\vdots &\vdots& \vdots \\
 0 & 0 &0 & \dots&1&\dots&\ast&|&\ast&\dots&\ast&\ast&\ast &\dots& 0 &0 \\
 \vdots & \vdots &\vdots&\vdots& \vdots&\dots & \vdots &\vdots&\vdots& \dots&\vdots & \vdots & \vdots & \vdots &\vdots& \vdots \\
0 & 0 &0 & \dots&0&\dots&\ast&|&\ast&\dots&\ast &\ast &\ast &\dots& \ast &0 \\
0 &0& 0&\dots &0&\dots&1&|&\ast&\dots&\ast &\ast & \ast &\dots &\ast& \ast  \\
 -&-&-&-&-&\dots&-&|&-&-&-&-&-&\dots&-&- \\
0 & 0&0 & \dots&0&\dots& 0&|&1&\dots&\ast &\ast & \ast &\dots &\ast & \ast \\
 \vdots & \vdots &\vdots&\vdots& \vdots&\dots & \vdots &\vdots&\vdots& \dots&\vdots & \vdots & \vdots & \vdots &\vdots& \vdots \\
0 & 0 &0 & \dots &0&\dots&0&|&0&\dots&1&\ast &\ast & \dots& \ast & \ast
\end{array}\right)}.
$$
The matrix $B_1$ has size $(m-1-t)\times (\Delta(\overline{G})+1+r),  0< r<t$\  where  $m-1-t=\Delta(\overline{G})+1+r$. This matrix  contains  the columns   of  the unknowns  $k_{m,1},\dots,k_{m,\Delta(\overline{G})+1},\\ (-g_{m,1}),\dots,(-g_{m,r})$.The matrix  $B_2$ has  size $(m-1-t)\times (\Delta(\overline{G})+1)$.  The $0$  matrix has  size  $t\times (\Delta(\overline{G})+1+r)$.  The matrix  $R$ has  size $t\times (\Delta(\overline{G})+1)$.

Since $t\le\Delta(\overline{G})<\Delta(\overline{G})+1$ the  system  $R\overrightarrow{g}=0$ has  infinitely many  solutions with  at  least  one  free  variable. By  the same  argument  as  in  case  2 we get  that the  system $S$ is consistent   and  the solution  with nonzero values  for $k_{m,1},\dots k_{m,\Delta(\overline{G})+1}$  gives  a vector $\overrightarrow{v}_m$ which  satisfies  all the adjacency conditions  and  orthogonal conditions.

\end{enumerate}
Hence, $Y_{m}$ has an  orthogonal  representation  of  vectors  in  $\Re^{\Delta(\overline{G})+1}$ satisfying   the conditions  (1) and (2). Thus $Y_{|G|}=G$ has an orthogonal representation  in $\Re^{\Delta(\overline{G})+1}$
satisfying  conditions (1) and (2).

Using  the  same argument in  the construction  of  vector $v_4$,  we    prove   that  $v_1,\dots,v_m$  is pairwise linearly  independent  set  of  vectors   in  $\Re^{\Delta(\overline{G})+1}$.
Finally  since $\msr(G)$ is  the smallest  dimension  in which  $G$ has an orthogonal  representation
$
\msr(G)\le \Delta(\overline{G})+1.
$
Since $\delta(G)+\Delta(\overline{G})= |G|-1$ we conclude that $G$ satisfies the delta conjecture, namely,
$
\msr(G)\le |G|- \delta(G).
$
\hfill $\Box$

\begin{observation}
In  the  construction  of  orthogonal representation of the  induced  graph $Y_m$  it is  sufficient  to consider   $t\le\left\lfloor\frac{m}{2}-1\right\rfloor $  if   $m$  is  even  and  $t<\left\lfloor\frac{m-1}{2}\right\rfloor $ if  $m$ is  odd.  In both  cases   we obtain  the  condition  $t< (m-1-t)$  that  we  need  to get  infinitely many  solutions    for  the  system $R\overrightarrow{g}=0$. This   difference  in  the upper  bounds   for  $t$  is important   for  small  values  of  $m$ but   for  larger  values  of $m$ these upper bounds are  asymptotically  equivalent. However,  it means that we  could  get  an orthogonal representation of  pairwise linearly independent  vectors  in $\Re^{\Delta(\overline{G})+1}$  for  some   graphs   which  are not  necessarily $\delta$-graphs.
\end{observation}
\begin{observation}
Reducing  the matrix $A$  to an echelon  form  needs  a finite  number of operations  as well as   reducing  the matrix $R$ to an  echelon form. It means   that all  the values  $k_{i,j}, i=1,2,\dots,m-1, j=1,2,\dots,\Delta(\overline{G})+1$ can be chosen   from different  field  extensions in such a way   that all the values  $\ast$ in the reduced echelon  form  of $A$    are nonzero and  belong to different field  extensions.
\end{observation}
\begin{observation}
The condition  of   choosing  values  $k_{i,j},  i=1,2,\dots,|G|, j=1,2,\dots,\\ \Delta(\overline{G})+1$ from  different  field  extensions was imposed  to  guarantee  the  consistency  of  the non-homogeneous system $S$. Also,  we use  this nonzero entries of  the vectors $\overrightarrow{v}_1,\dots, \overrightarrow{v}_{m-1}$  to guarantee  the adjacency conditions and  orthogonal conditions  of the  vector $\overrightarrow{v}_m$ corresponding to the  newly  added vertex. But calculating  the orthogonal representation  using  this approach  could  be time consuming. Since   we  know  that  it is possible  to get an orthogonal  representation  of  $\delta$-graph $G$   in  $\Re^{\Delta(\overline{G})+1}$  and since the  representation  is  not   unique,  it may be   possible  to  calculate   the orthogonal representation using  integers  or   rational numbers. However,   calculating  the orthogonal  representation of  a $\delta$-graph $G$ in  this  way  could  also be tedious because   we may need  to apply a backtracking  procedure during the calculation   due  to  some   adjacency  conditions   of  the  vector corresponding  to  the  newly   added  vertex may not  be satisfied. When  that  happens,  we may  need to go back to  some of the  previous  vectors  and  recalculate  them   until    we  fix   the adjacency conditions.
\end{observation}
\section{Examples  of  $\delta$-graphs  and their  $\msr$}
The  result  proved above   give us a  huge family  of  graph  which  satisfies  delta conjecture. Since,  the   complement    of   a  C-$delta$ graphs  is  a $\delta$-graph,  it  is  enough  to identify  a  C-$\delta$-graph    and  therefore  we know   that   its  complement  is a $\delta$-graph satisfying  delta  conjecture.

Some  examples   of  C-$\delta$ graphs that  we can  find  in  \cite{PD} are  the  Cartesian Product $K_n\square P_m,n\ge  3,  m\ge  4$, Mobi\"us Lader $ML_{2n}, n\ge 3$, Supertriangles $Tn, n\ge 4$,  Coronas $S_n\circ P_m, n\ge  2 , m\ge 1$ where  $S_n$  is a  star  and  $P_m$ a path, Cages like Tutte's (3,8) cage,  Headwood's (3,6) cage   and many others, Blanusa Snarks  of type  $1$ and  $2$  with  $26, 34$, and $42$ v\'ertices,  Generalized Petersen Graphs $Gp1$ to  $Gp16$, and  many  others.

In  order   to  show  the  technique  used  in  the proved  result   consider  the following  example
\begin{example}\label{upper2}

If  $G$ is  the Robertson's (4,5)-cage on 19 vertices  then  it is  a 4-regular C-$\delta$ graph. Since $\Delta(G)=4$, the $\msr(G)\le 5$.  To  see  this is a C-$\delta$ graph it  is enough  to label its  vertices  in the  way  shown  in the next  figure:
\begin{center}
\includegraphics[height=50mm]{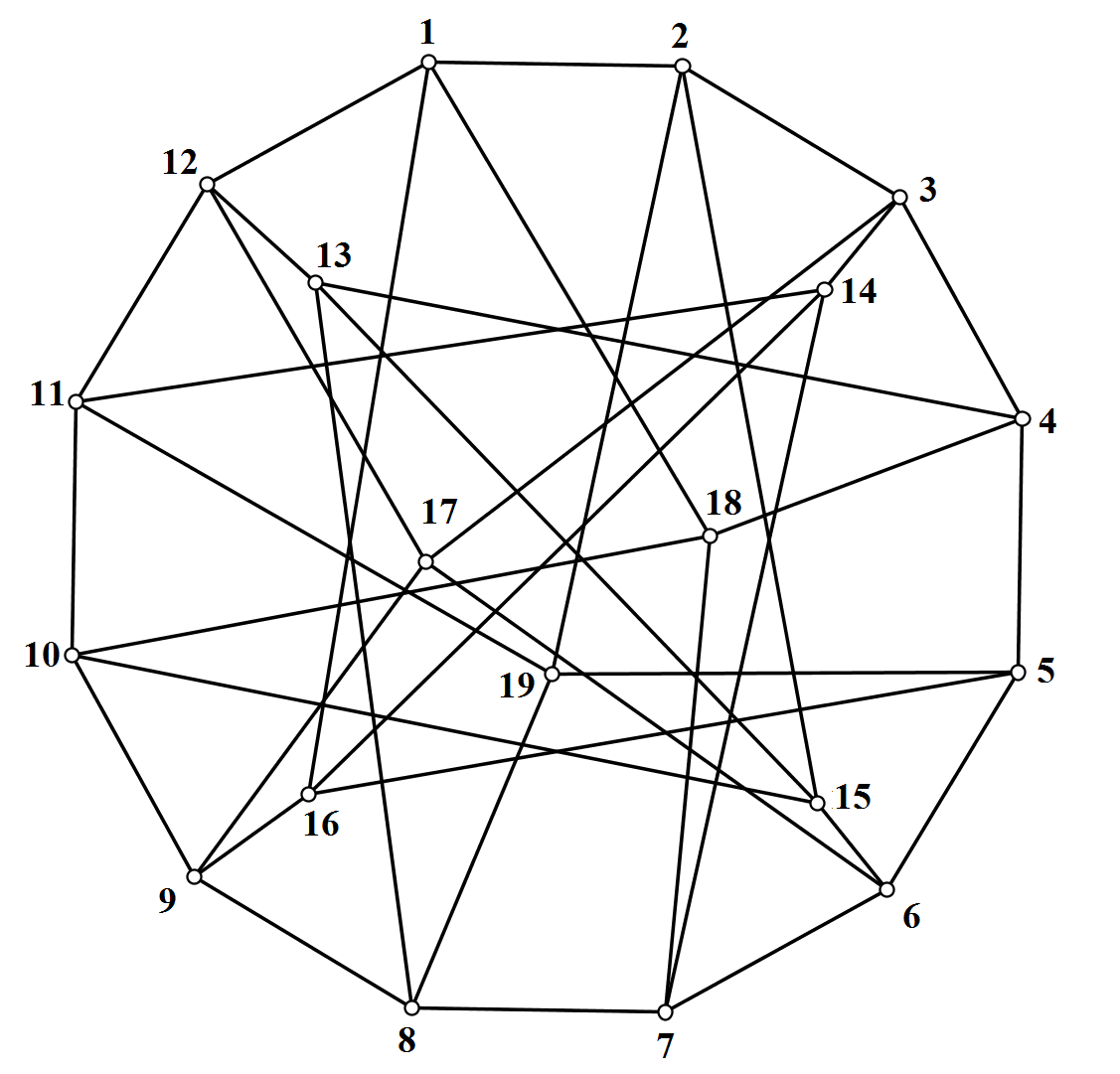}
 \end{center}
\vspace*{-0.1in}\begin{figure}[h]
\centering
%\caption
{Figure B.2 Robertson's (4,5)-cage (19 vertices)}
\label{figA.1.2}
\end{figure}
\end{example}
\section{Conclusion}

The  result  proved above   give us  a a tool  to  identify a  wide  range of  families  of  graphs  which  satisfy  $\delta$ conjecture.  The  techniques  used  in  the  proof   could  be  used   in  future  research   as  a new approach  to solve  delta  conjecture. However, it  is  clear  that  the main problem  is  still open.
\section{Acknowledment}
I  would like  to  thanks   to  my  advisor  Dr. Sivaram Narayan  for   his  guidance and  suggestions  of  this  research.  Also  I  want  to  thank  to  the math  department  of  University  of  Costa  Rica and  Universidad  Nacional Estatal  a Distancia  because  their  sponsorship  during my  dissertation  research  and    specially  thanks  to   the math  department    of  Central Michigan University   where  I  did  the  researh   for  this  paper.
%****************************************************************************
%*                        BIBLIOGRAPHY                                      *
%****************************************************************************
\renewcommand{\baselinestretch}{1} 

\end{document}